\documentclass[]{ametsoc}

\usepackage{float}
\usepackage{subfigure}
\newcommand{\bs}{\boldsymbol}

\journal{mwr}

\bibpunct{(}{)}{;}{a}{}{,}

\title{A data-driven method for improving the correlation estimation in serial ensemble Kalman filters}

\authors{Mich\`ele De La Chevroti\`ere \correspondingauthor{Department of Mathematics, the Pennsylvania State University, 109 McAllister Building, University Park, PA 16802-6400, USA}}
\affiliation{Department of Mathematics, the Pennsylvania State University, USA} 

\email{mdelachev@psu.edu}

\extraauthor{John Harlim}
\extraaffil{Department of Mathematics and Department of Meteorology and Atmospheric Science, the Pennsylvania State University, USA}

\abstract{A data-driven method for improving the correlation estimation in serial ensemble Kalman filters is introduced. The method finds a linear map that transforms, at each assimilation cycle, the poorly estimated sample correlation into an improved correlation. This map is obtained from an offline training procedure without any tuning as the solution of a linear regression problem that uses appropriate sample correlation statistics obtained from historical data assimilation outputs.
In an idealized OSSE with the Lorenz-96 model and for a range of linear and nonlinear observation models, the proposed scheme improves the filter estimates, especially when the ensemble size is small relative to the dimension of the state space.}

\begin{document}

\maketitle

\section{Introduction}

It is a well known fact  \citep{evensen:94} that when small ensemble sizes are used in an ensemble Kalman filter (EnKF), the estimated correlation tends to overestimate the correlations between two variables that are far away even if the true correlation is small \citep{lorenc:03}. The typical approach to address this issue is by zeroing (and damping) the correlation coefficients of longer (intermediate) distances, which practically means ignoring (weighting less) such observations. Such practice is known as \textit{localization} \citep{houtekamer:98}. Many localization functions have been proposed, with tapering \citep{furrerbengtsson:07}, the Gaspari-Cohn function \citep{gaspari:99}, or other distance-based localization functions for multiphase flows \citep{chenoliver:10}. While these approaches have been very useful, they require exhaustive tuning which can be expensive for problems with multiscale spatial correlations \citep{dong:11,ki:14,flowerdew:15}. Furthermore, since the correlation between a state variable $\boldsymbol{x}$ and a nonlinearly observed state $h(\boldsymbol{x})$ may not necessarily be a function of spatial distance, specifying a localization function for such correlations may not be trivial. Alternatively, a class of adaptive (or ``flow-dependent") methods has also been proposed, however they depend on many parameters such as the integer powers in  \citet{bishophodyss:07}, or require specification of a prior distribution for the sampling correlation \citep[e.g. see][]{zz:14}. There is also a nonparametric method \citep{jun:11} which provides a guideline for estimating the background covariance matrix.

More recently, an approach for estimating the localization function that requires almost no tuning has been proposed by \citet{andersonlei:13}. This approach consists of training empirical localization functions (ELF) from an observation system simulation experiment (OSSE) data assimilation output by minimizing the analysis error, i.e., the difference between the analysis mean update and the truth (or observations when the truth is unknown). Numerically, they obtain ELF in stages: The first ELF is obtained by minimizing the analysis error from an EnKF without localization, the second ELF is obtained by minimizing the analysis error from EnKF using the first ELF estimate, and so on. In this sense, there is still some tuning required, for instance in determining the stopping criteria. 

Motivated by this approach, we introduce a data-driven method for improving the sample correlation estimation in the serial EnKF when small ensemble sizes are used. In particular, we consider finding a linear map that takes poorly estimated sample correlations from each EnKF cycle and transforms them into an improved correlation matrix. This linear map is computed by regressing appropriate sample correlation statistics from historical data assimilation output obtained in a training phase and we will show that there is no tuning required. We shall see that this linear map generalizes the standard localization function, and our method for specifying this map can be interpreted as a discrete (or Monte-Carlo) approximation to a minimization problem with cost function inspired by \citet{anderson:12} and \citet{flowerdew:15}. Our approach is relatively easy to implement in any observation network and it does not parameterize the prior distribution as in \cite{zz:14} or the likelihood function as proposed in \citet{anderson2016}. 

The remainder of this paper is organized as follows: In Section~\ref{section2}, we review the basic statistical results on sampling error in estimating correlation which motivate the current approach. In Section~\ref{section3}, we introduce our new algorithm for specifying the linear map and describe its implementation with a serial EnKF of \citet{anderson:03}. In Section~\ref{section4}, we demonstrate the numerical results in an idealized OSSE setup using the Lorenz-96 model \citep{lorenz1996predictability}, and compare it to an EnKF with the well-tuned Gaspari-Cohn localization function as well as ETKF with a very large ensemble size, where we use the ETKF formulation of \cite{hunt:07} without localization.  We close this paper with a short summary in Section~\ref{section5}.

\section{Sampling error in correlation statistics of EnKF}\label{section2}

The goal of this section is to give readers a brief review to understand sampling errors in an EnKF. In particular, we will show that it is rather difficult to provide a practically useful estimate for sampling errors in the EnKF in the linear and Gaussian case, and discuss how this issue is compounded in high-dimensional nonlinear filtering problems. While the detailed formulation in this discussion may seem to be a digression, it aims to provide intuition for the algorithm proposed in Section~\ref{section3}.

Consider the discrete-time linear and Gaussian filtering problem, 
\begin{eqnarray}
\bs{x}_{m+1} &=& F \bs{x}_m + \bs{\zeta}_{m+1}, \quad \bs{\zeta}_m  \sim\mathcal{N}(0,Q),\label{lineardynamics} \\
\bs{y}^o_m &=& H\bs{x}_m + \bs{\epsilon}_m, \quad \bs{\epsilon}_m\sim\mathcal{N}(0,R),\label{linearobs}
\end{eqnarray}
where $F$ denotes the deterministic model operator and $Q$ denotes the Gaussian system noise covariance matrix. The time index is defined as $t_{m+1}=t_m+\tau$, where $\tau$ denotes the observation time step. 
The optimal solution (in the least-square sense) to this filtering problem, assuming both noises are uncorrelated, $\mathbb{E}[\bs{\epsilon_m\zeta_n}^\top]=0$ for all $m,n$, is given by the Kalman filter formula,
\begin{eqnarray}
\overline{\bs{x}}^{\,a}_m &=& \overline{\bs{x}}^{\,f}_m + \mathcal{K}_m(\bs{y}^o_m-H\,\overline{\bs{x}}^f_m), \quad\quad P^{\,a}_{m} = (I-\mathcal{K}_mH)P^{\,f}_m, \nonumber\\
\mathcal{K}_m &=& P_m^{\,f}H^\top(HP_m^{\,f}H^\top+R)^{-1}\label{kalmanfilter}\\
\overline{\bs{x}}^{\,f}_{m+1} &=& F \overline{\bs{x}}^{\,a}_m, \quad\quad P^{\,f}_{m+1} = FP^{\,a}_{m}F^\top + Q  \nonumber
\end{eqnarray}

where $\overline{\bs{x}}^{\,f}_m$, $P^{\,f}_m$ denote the forecast mean and covariance estimates before accounting for the observations $\bs{y}^o_m$, and $\overline{\bs{x}}^{\,a}_m$, $P^{\,a}_m$ denote the analysis mean and covariance estimates after accounting for the observations $\bs{y}^o_m$. In this rather restrictive scenario, the time asymptotic limit of the Kalman filter covariance solutions, $P^{\,a}_m$, is the solution $P^*$ of the  discrete Ricatti equation
\begin{equation}
P^* = (I-\mathcal{K}H)(FP^* F^\top+Q),\nonumber
\end{equation}
where ${\cal K}=P^*H^\top(HP^*H^\top+R)^{-1}$ denotes the limiting Kalman gain, i.e. we have $P^{\,a}_m \to P^*$ as $m\to\infty$.

The main idea behind ensemble Kalman filters \citep[e.g.][]{evensen:94,anderson:01,bishop:01} is to approximate the forecast distribution by an ensemble of solutions, use the Kalman filter formula in \eqref{kalmanfilter} to transform the first two empirical prior moments, $\hat{\overline{\bs{x}}}^{\,f}_m$, $\hat{P}^{\,f}_m$ (obtained from the ensemble average) to analysis empirical mean and covariance, $\hat{\overline{\bs{x}}}^{\,a}_m$, $\hat{P}^{\,a}_m$, and finally draw a Gaussian sample $\{\bs{x}^{a,k}_m\}$ such that
\begin{eqnarray}
\hat{\overline{\bs{x}}}^{\,a}_m = \frac{1}{K}\sum_{k=1}^K \bs{x}^{a,k}_m\quad \mathrm{and} \quad
\hat{P}^{\,a}_m = \frac{1}{K-1}{X}^a_{m}\, {{X}^a_{m}}^\top,\nonumber
\end{eqnarray}
where each column of the matrix ${X}^a_{m}$ consists of the perturbation of the $k^{\mathrm{th}}$ ensemble member from it's sample mean, $\bs{X}^{a,k}_{m} = \bs{x}^{a,k}_m-\hat{\overline{\bs{x}}}^{\,a}_m$, and $K$ denotes the ensemble size. Here, $m$ is the time index.

In the Gaussian and linear case, the EnKF solutions are optimal when the sample mean estimate, $\hat{\overline{\bs{x}}}^{\,a}_m$, matches the Kalman mean estimate $\overline{\bs{x}}^{\,a}_m$ and the sample covariance estimate, $\hat{P}^{\,a}_m$, matches the asymptotic covariance $P^*$ as $m\to\infty$. In such a case, the EnKF analysis ensemble perturbations, $\bs{X}^{a,k}_{m}$, are samples of a time independent Gaussian distribution with mean zero and covariance matrix $P^*$ or correlation matrix $\rho = D^{-1/2}P^* D^{-1/2}$, where $D$ is a diagonal matrix whose diagonal components are the diagonal components of $P^*$.  
That is, as $m\rightarrow \infty$, the EnKF numerically estimates $P^*$ using the sample covariance matrix $\hat{P}^a_m$ at each data assimilation step $m$.

From basic statistical theory \citep{cb:2002}, we know that the standard deviation of the estimated variance $\hat{P}^{\,a}_{ii}$ obtained from $K$ independent samples of a scalar Gaussian random variable $x_i$, the $i^\mathrm{th}$ component of $\bs{x}$, is $\sigma_i = \sqrt{\frac{2}{K-1}} P^*_{ii}$, which is smaller than the variance $P^*_{ii}$ when the ensemble size $K>3$. In operational implementations, where the ensemble size $K\approx 10-100$, this standard deviation is quite small. Therefore, it is reasonable to assume that the sample covariance $\hat{P}^{\,a}_{ii}$ is an accurate estimation of the true variance $P^*_{ii}$. A similar assumption is taken in  \citet{anderson:12, anderson2016}. 

Now let's turn our attention to the sampling errors in the estimation of the non-diagonal components of $P^*$ or essentially the correlation matrix $\rho$, assuming that the diagonal component of $P^*$ can be accurately estimated by the EnKF solutions $\{\bs{x}^{a,k}_m\}_{k=1,\ldots, K}$. In other words, we are estimating $\rho$ with a sample correlation matrix $r^K_m \equiv D^{-1/2} \hat P^a_m D^{-1/2}$, assuming that $D=\mathrm{diag}(P^*)=\mathrm{diag}(\hat P^{\,a}_m)$, using the EnKF with an ensemble size $K$. Of course these samples are not temporally independent for observation time interval that satisfies $\tau < T_c$,  where $T_c$ denotes the decaying time scale of the underlying dynamics. For $m\to\infty$ and large enough $\tau>T_c$ such that $\bs{X}^{a,k}_{m}$ are temporally independent samples of a multivariate Gaussian distribution with mean zero and covariance $P^*$, the distribution and moments of its sample correlation coefficient $r^K_m$ are explicitly known \citep{hotelling:53}. In particular, $\mathbb{E}[(r_{{\,ij}}^K-\rho_{ij})^2] = (1-\rho_{ij}^2)^2K^{-1}+\mathcal{O}(K^{-2})$, which means that the error is on the order of $K^{-1}$ when $\rho_{ij}=1$ and $K^{-1/2}$ otherwise. Of course, this error estimate is valid only when $\bs{x}^{a,k}$ are independent samples of a Gaussian prior and the filtering problem is linear, and so it is not useful in practice.

In general high-dimensional non-Gaussian and nonlinear settings, the problem becomes much harder since the true correlation $\rho_m$ is time-dependent (it does not equilibrate to $\rho$). Furthermore, the true correlations $\rho_m$ are statistics with respect to the true filtering distribution $p(\bs{x}_m|\bs{y}^o_m)$ which is almost always non-Gaussian in the general nonlinear filtering problem. The bottom line is that if one approximates the nonlinear filtering problem with an EnKF, the sample correlation is not the same as the true correlation even in the limit of large ensemble size,
\begin{equation}
\lim_{K\to\infty} r^K_m \equiv r_m \neq \rho_m, \quad\forall m,\nonumber
\end{equation}
where $r^K_m$ denotes the sample correlation obtained from an EnKF with ensemble size $K$, $r_m$ denotes the limiting correlation as $K\to\infty$, and $\rho_m$ denotes the true correlation. While estimating the true non-Gaussian statistics $\rho_m$ may be very difficult in practice, it is not hopeless to estimate the limiting sample correlation $r_m$ from a sample correlation $r^K_m$ that is obtained from an EnKF with a finite ensemble size $K$. This is precisely what we will discuss in the next section, that is, we propose to find a linear map $\mathcal{L}$ that takes the sample correlation $r^K_m$ to the limiting sample correlation $r_m$. In a special case, this map is exactly what is commonly known as the \textit{localization function}.

\section{Algorithm}\label{section3}

The main idea in \citet{anderson:12} and \citet{flowerdew:15} is to find a localization function by minimizing the mean square difference of the analysis from the truth. Motivated by their ideas and the theoretical understanding above, we propose to find a linear map $\mathcal{L}\in\mathbb{R}^{N\times N\times M}$, a three-dimensional array which transforms the sample correlations between $\bs{x}$ and $\bs{y}$, to the corresponding limiting sample correlation. We will learn this map from historical data assimilation outputs as discussed in Section 3\ref{section31} below. Later on,  in Section 3\ref{section32}, we show a straightforward modification of a serial EnKF \citep{anderson:03} using this linear map $\mathcal{L}$.  

\subsection{Learning the ``localization" map}\label{section31}

First, let us set our notation. To simplify the notation in the remainder of this paper, we suppress the subscript time $m$ and only use it when necessary. Since we are interested in filtering problems with general observation functions $h:\mathbb{R}^N\rightarrow\mathbb{R}^M$, we consider a sample cross correlation between $\bs{x}$ and $\bs{y}$ as follows,
\begin{equation}
r_{\boldsymbol{x}\bs{y}}^K = D^{-1/2}_{X^a} \Big(\frac{1}{K-1}X^a{Y^a}^\top\Big) D^{-1/2}_{Y^a} \in\mathbb{R}^{N\times M}, \label{crosscorrelation}
\end{equation}
where $X^a$ is defined as before, and the columns of $Y^a$ are defined as $\bs{Y}^{a,k} = h(\bs{x}^{a,k})-\hat{\overline{\bs{y}}}^{\,a}$, where $\hat{\overline{\bs{y}}}^{\,a}=\sum_{k=1}^K h({\bs{x}}^{a,k})/K$. In \eqref{crosscorrelation}, the matrices $D_{X^a} = \mathrm{diag}(X^a{X^a}^\top)/(K-1)$ and $D_{Y^a} = \mathrm{diag}(Y^a{Y^a}^\top)/(K-1)$ are diagonal. Notice that in the case of linear observation operators $h(\bs{x})=H\bs{x}$, we have $Y^a = HX^a$.  In the remainder of this article, we suppress the subscript $\bs{x}\bs{y}$ and understand that {\color{black}$r^K = r_{\bs{x}\bs{y}}^K$} with limiting correlation $r=\lim_{K\to\infty}r^K$. We also use the notation $\mathcal{L}(\cdot,i,j) \in \mathbb{R}^N$ to denote an $N$-dimensional column vector corresponding to fixing the second and third component in the array $\mathcal{L}$ by $i$ and $j$, respectively.

We consider finding a linear functional $\mathcal{L}(\cdot,i,j):[-1,1]^N\to[-1,1]$ that takes the poorly sampled correlations between all state variables and the observation $j$, and maps them to an improved correlation $\tilde{r}(i,j)$ for the pair of state variable $i$ and observation $j$, defined as $\tilde{r}(i,j) = {r^K(\cdot,j)}^\top \mathcal{L}(\cdot,i,j)$.
We propose to find this map by solving
\begin{equation}
\min_{\mathcal{L}(\cdot,i,j)} \int_{[-1,1] \times [-1,1]} ({r^K(\cdot,j)}^\top \mathcal{L}(\cdot,i,j)-r(i,j))^2 p(r^K|r)p(r) \,dr^K\,dr,\label{optimization}
\end{equation}
such that $\tilde{r}(i,j)$ best approximates $r(i,j)$, the limiting correlation between the state variable $i$ and the observation $j$. Here $p(r)$ is the distribution of the limiting sample correlation whereas $p(r^K|r)$ is the conditional distribution of the sample correlation $r^K$ given $r$, where we assume that these distributions are stationary. In this formulation, the goal is to find an estimator for the limiting correlation $r(i,j)$ through a linear map which combines information from all the sample correlations $r^K(\cdot,j)$, whereas the approach in \citep{anderson:12,flowerdew:15}  considers a scalar factor applied to the corresponding single  sample correlation $r^K(i,j)$ only. Their approach is similar to assuming that the map $\mathcal{L}(\cdot,\cdot,j)\in\mathbb{R}^{N\times N}$ is diagonal in our formulation (see the discussion at the end of this section).

To realize {\color{black}Eqn}.~\eqref{optimization} numerically, we consider the following training procedure. Consider an OSSE historical data assimilation analysis ensemble, $\{\bs{x}^{a,k}_m\}$, where the index $m=1,\dots,T$ denotes time and the index $k=1,\dots,L$ denotes ensemble member, obtained from a large ensemble of size $L\gg 1$, such that $r^L\approx r$. Here, $T$ denotes the number of data assimilation cycles in the training period. Given this analysis ensemble, one can sample the limiting correlation by $r_m^L$, and also obtain subsampled correlations $\{r_{m,s}^K \}$, using only $K$ members out of $L$. Here we use the index $s=1,\dots,S$ to denote the number of possible ways to choose $K$ members out of $L$, thus one can sample as many as $\mathrm{max}(S) = {L \choose K}$. Therefore, from the historical analysis ensemble, we obtain samples $r^L_m \sim p(r)$ and  $r^K_{m,s} \sim p(r^K|r=r^L_m)$. We can then apply a Monte-Carlo approximation to the integral in \eqref{optimization} such that the minimization problem becomes
\begin{equation}
\min_{\mathcal{L}(\cdot,i,j)} \frac{1}{TS}\sum_{m,s=1}^{T, S} ({r^K_{m,s}(\cdot,j)}^\top \mathcal{L}(\cdot,i,j)-r^L_m(i,j))^2,\label{regression}
\end{equation}
which is essentially a least square problem.  To see this, define $A\in\mathbb{R}^{TS\times N}$, where each row of $A$ is $r^K_{m,s}(\cdot,j){^\top}$, $\bs{u} = \mathcal{L}(\cdot,i,j)\in\mathbb{R}^N$, and the vector $\bs{b}\in\mathbb{R}^{TS}$ whose components are $r^L_{m,s}(i,j)$, where $r^L_{m,s}(i,j) = r^L_{m}(i,j)$ for any $m,s$. Then the linear regression problem in \eqref{regression} can be rewritten in compact form as $\min_{\bs{u}} (A\bs{u}-\bs{b})^\top(A\bs{u}-\bs{b})$, ignoring the scalar constant factor $(TS)^{-1}$. An explicit solution to this linear regression problem is given by $\bs{u}=(A^\top A)^{-1}A^\top \bs{b}$. Repeating this procedure for every indices $i$ and $j$, we obtain the linear map $\mathcal{L}$.

In this formulation, there are parameters $S, K, T, L$. For large enough $T$, we can indeed set $S=1$ and still obtain accurate filter estimates as we will show below. The choice of the parameter $K$ depends on the ensemble size that will be used in the actual online data assimilation. For the parameters $T$ and $L$, one can choose $L$ as large as possible depending on the computing capacity and $T$ as large as possible depending on the length of the training data set or historical correlation statistics. In contrast with \citet{andersonlei:13}, we do not need the historical true state to specify the map $\mathcal{L}$. In this sense, this data-driven method requires no tuning.

The complexity in constructing the correlation $r_{\boldsymbol{x}y_i}$ between a single scalar observation $y_i$ and the state variable $\boldsymbol{x}\in\mathbb{R}^N$ using an ensemble of size $L$ is $\mathcal{O}(LN)$. Similarly, the cost for constructing the subsampled correlation using $K$ members, where $K\leq L$, is $\mathcal{O}(KN)$. If these calculations are performed on historical data of length $T$ with $S$-subsamples, then the complexity in constructing the regression matrix $A$ is $\mathcal{O}(LNTS)$. In each minimization problem, the complexity in forming $A^\top A$ is $\mathcal{O}(N^2TS)$, in forming $A^\top b$ is $\mathcal{O}(NTS)$, and in inverting $A^\top A$ (say with LU factorization) is $\mathcal{O}(N^3)$. To avoid an underdetermined linear problem, $TS>N$, so $\mathcal{O}(LNTS)$ dominates, assuming that $L>N$. In our numerical experiment below, we will show accurate results even for $S=1$. For high-dimensional problems one may consider a local least square fit, by assuming that each observation is only correlated to its nearby state space of dimension $N_{loc}<N$, as opposed to the whole state space. With these considerations, the complexity in solving each regression problem reduces to $\mathcal{O}(LN_{loc}T)$ and one can choose much smaller $L$ to satisfy $N_{loc}<L<N$.

In the special case where $\mathcal{L}(\cdot,\cdot,j)\in\mathbb{R}^{N\times N}$ is restricted to be diagonal, the proposed minimization problem \eqref{optimization} seeks to estimate the localization function such that $\mathcal{L}(i,i,j) r^K(i,j) \approx r(i,j)$, where $\mathcal{L}(i,i,j) \in \mathbb{R}$; we refer to the resulting matrix from fitting to the diagonal $\mathcal{L}$ as $\mathcal{L}_d$. In this case, the minimization procedure is pointwise and corresponds exactly to the strategies proposed in \citep{anderson:12,andersonlei:13,flowerdew:15}, with the only difference being in how the densities are estimated from the training data.

\subsection{A serial EnKF with transformed correlation functions}\label{section32}
In this paper, we will consider the implementation of the map $\mathcal{L}$ on the serial local-least-squares LLS-EnKF \citep{anderson:03}, which assumes a diagonal observation error covariance $R$ and uses the unperturbed observation component $y^o_j$ to sequentially update the ensemble solutions. In each assimilation step, the LLS-EnKF executes the following update
\begin{eqnarray}
 \hat{\overline{y}}^{\,a}_j &=& \hat{\overline{y}}^{\,f}_j + (\hat{P}^f_{y_jy_j}+R_{jj})^{-1}\hat{P}^f_{y_jy_j} (y^o_j-\hat{\overline{y}}^{\,f}_j),\nonumber\\
y^{a,k}_j &=& \hat{\overline{y}}^{\,a}_j + \sqrt{\frac{R_{jj}}{R_{jj}+\hat{P}^f_{y_jy_j}}} (y^{f,k}_j-\hat{\overline{y}}^{\,f}_{j})\nonumber
\end{eqnarray}
and regresses the update $\Delta y^k_j = y^{a,k}_j -y^{f,k}_j $ onto the state variables as follows,
\begin{equation}
\boldsymbol{x}^{a,k} = \boldsymbol{x}^{f,k} + \frac{\hat P^{\,f}_{\bs{x}y_j}}{{\hat P^{\,f}}_{y_jy_j}}\Delta y^k_j.\label{lls} 
\end{equation}
The straightforward modification on the LLS-EnKF is to replace
\begin{equation}
\hat{P}^f_{\bs{x}y_j} = D_{X^f}^{1/2} r^K(\cdot,j) \sqrt{{\hat P^f}_{y_jy_j}} \longleftarrow   D_{X^f}^{1/2} \mathcal{L}(\cdot,\cdot,j)^\top r^K(\cdot,j) \sqrt{{\hat P^f}_{y_jy_j}},\label{corr_transfo}
\end{equation}
in \eqref{lls}, where $D_{X^f}=\mathrm{diag}(X^f{X^f}^\top)/(K-1)$ and $r^K(\cdot,j)$ are estimated from the EnKF with  ensemble size $K$.

\section{Numerical experiments on Lorenz-96 model}\label{section4}

We test our methodology on an OSSE with the Lorenz-96 model \citep{lorenz1996predictability}
\begin{equation}
\dfrac{dx_j}{dt} = (x_{j+1} - x_{j-2}) x_{j-1} - x_j + \mathcal{F}, \quad j = 0, \dots, N-1, \label{lorenz} 
\end{equation}
with $N=40$ state variables ($\boldsymbol{x} \in \mathbb{R}^{40}$) equally spaced on a periodic one-dimensional domain of length 40 ($x_{N} = x_0$ and $x_{N+1} = x_1$). The forcing parameter $\mathcal{F}=8$ is set  for a strongly chaotic regime. A realization of the ``true'' signal or ``nature run'' $\boldsymbol{x}^\dagger$ is obtained by integrating the system (\ref{lorenz}) forward in time, from random initial conditions, using a fourth-order Runge-Kutta scheme with a time step $\Delta t = 0.05$, or 6 hours. 

Synthetic observations are generated from that truth trajectory  based on the model
\begin{equation}
\boldsymbol{y}^o_m=  h\left(\boldsymbol{x}^\dagger_m \right) + \boldsymbol{\epsilon}_m, \qquad \boldsymbol{\epsilon}_m  \sim \mathcal{N}(0,I), \label{obs_model}
\end{equation}
at discrete time $t_m$, where $\boldsymbol{\epsilon}_m$, the observation noise, is unbiased and Gaussian with identity covariance matrix. Observations are made at uniform discrete times $t_{m+1}-t_m = n \Delta t$, where $n\in \mathbb{N}$ is the observation time step. In total 30 000 observations are collected, with the first 10 000 used for training purposes, and the last 20 000 for diagnostic evaluation. We will test both linear and nonlinear observation operators $h: \mathbb{R}^N \rightarrow \mathbb{R}^M$, and in particular will focus on the following observation configurations:
\begin{itemize}
\item \textit{Linear direct observations}: every $p^{\mathrm{th}}$ grid point is observed, where $p=N/M$, $M\in \{10,20,40 \}$. 
\item \textit{Linear indirect observations}: every other grid point is observed indirectly as a sum of itself and the 6 nearest neighbours. In this case, the linear observation operator $H$ is defined such that 
\begin{equation}
(H\boldsymbol{x})_{j} = \sum_{k=-3}^{3} x_{mod(2j+k,N)}, \quad j = 1,\dots,M \quad (M=20),\nonumber
\end{equation}
where the modulus is to preserve the periodicity of the domain.

\item \textit{Nonlinear indirect observations}: We consider a nonlinear observation function of the form
\begin{equation}
h_j(\boldsymbol{x}) =  \sum_{k=-3}^3 w(x_{j+k}) x_{j+k},\qquad j = 1, \dots, M \quad (M=10),\label{nlinobs}
\end{equation}
with a state dependent weight, $w: [a,b] \rightarrow [0,1]$ defined as
\begin{equation*}
w(x_{j+k}) = \frac{c_{k+4}}{2}\Big[1+\cos\left(  \dfrac{2\pi}{b-a} \Big(   \dfrac{}{}x_j -\dfrac{a+b}{2}  \Big)\right)\Big], 
\end{equation*}

where $\boldsymbol{c} = (1,0.8,0.4,0,0.4,0.8,1)$ so when $k=\pm 3$, $c_{k+4}=1$, etc. Here, the parameters $a$ and $b$ are the extrema of the data $\{\boldsymbol{x}_m\}_{m=1,\ldots,30 000}$.
\end{itemize}

\subsection{Localization map $\mathcal{L}$}

Given some estimates $r$ and $r^K$, the method in Section 3 provides an algorithm to learn the map $\mathcal{L}$ offline with little computational overhead, by solving an appropriate minimization problem. Indeed for every pair $(i,j)$ of and state variable $i$ and observation $j$, the map $\mathcal{L}(\cdot,i,j)$ is given as the solution of the least square problem (\ref{regression}), which is effectively a linear regression of the sample background correlation $r^K$ onto $r$. 

\begin{figure}
\centerline{
\hspace{0.3cm}\includegraphics[width=4cm,height=4.5cm]{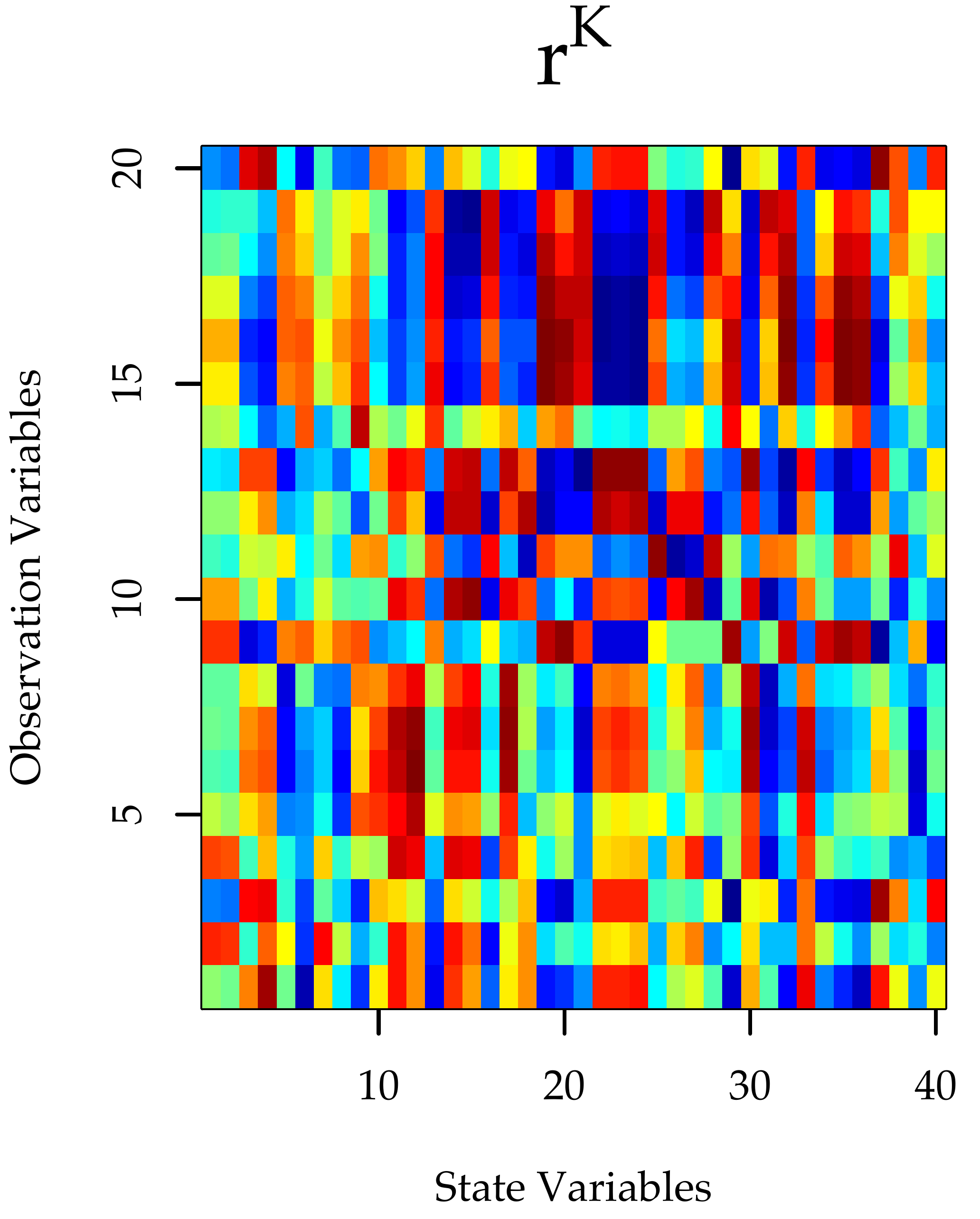}\includegraphics[width=4cm, height=4.5cm]{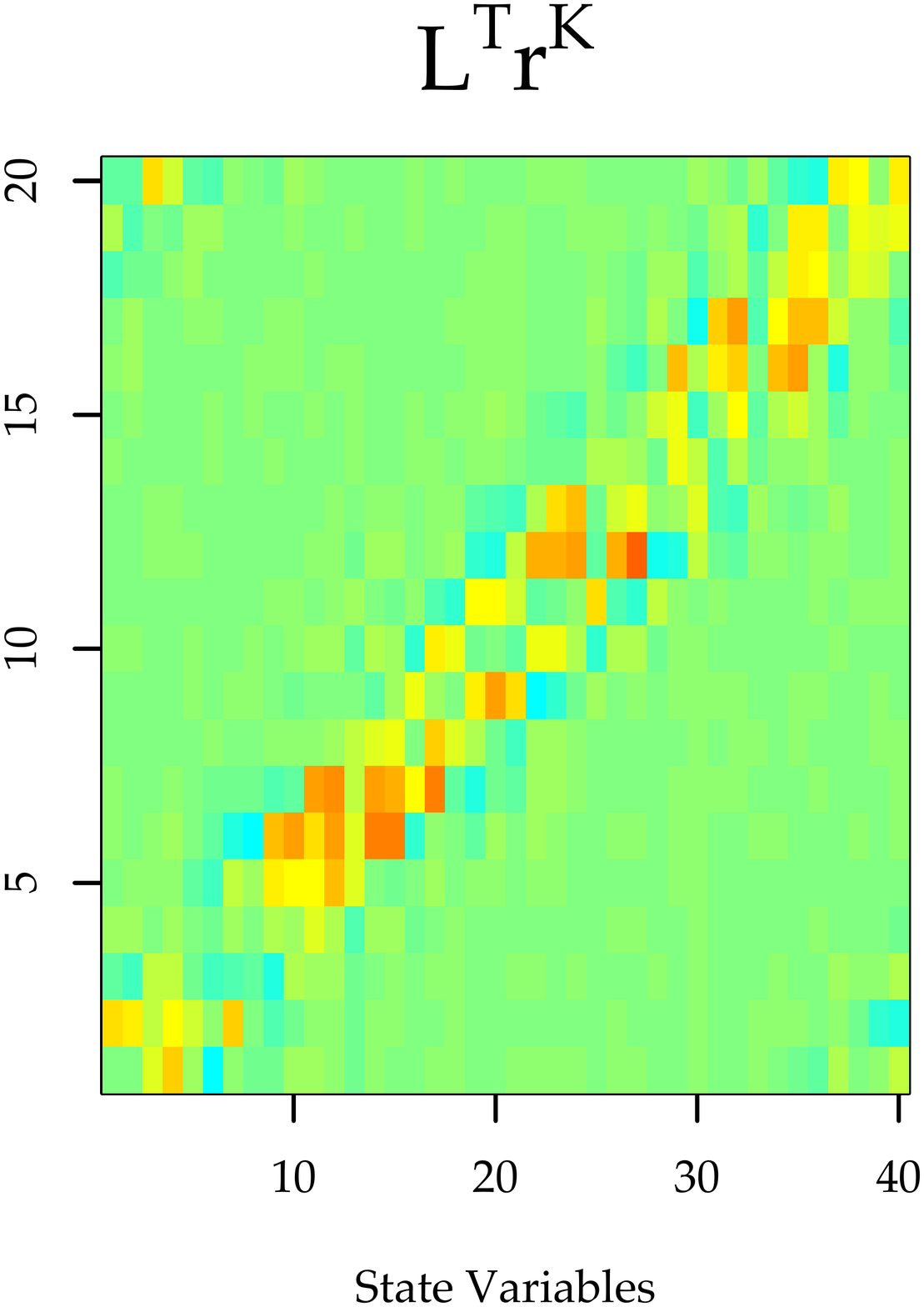}
\includegraphics[width=4.1cm, height=4.5cm]{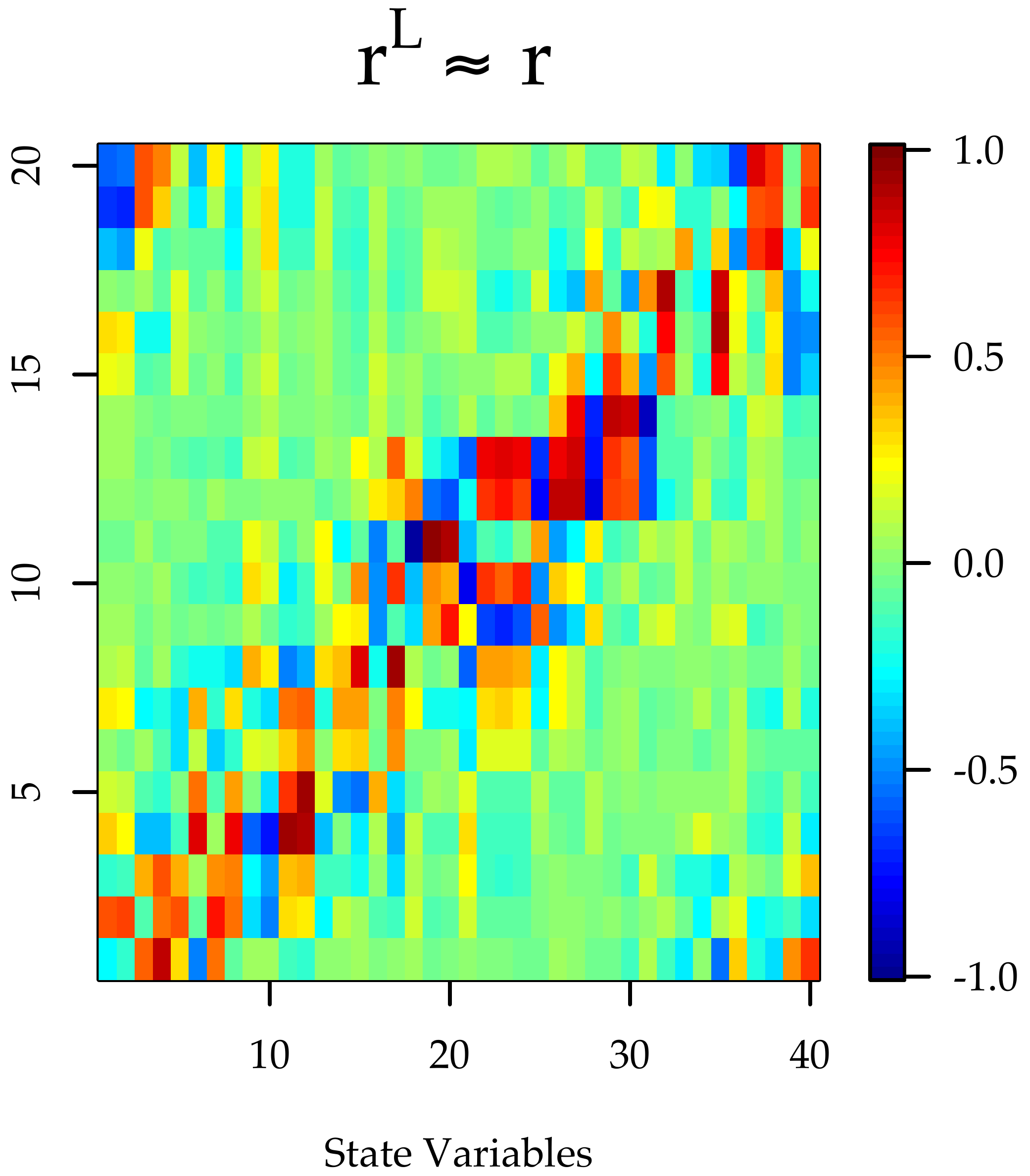}
}

\caption{The map $\mathcal{L}$ transforms the undersampled background error correlation matrix $r^K$ into a sample of correlation matrix  $\mathcal{L}^T r^K$ resembling more the regressor $r^L$. A snapshot of all three correlation matrices are shown here for 20  linear direct observations taken every model time step, for an ensemble of size $K=5$. The regression is based on the output sample correlation of an ETKF with ensemble size $L$=500.} \label{old_new_corr}

\end{figure}

For the OSSE data assimilation system described by the model equations (\ref{lorenz})-(\ref{obs_model}), the limiting correlation $r$ needs to be estimated in some way. Here we use the output sample correlation of an ETKF with ensemble size 500 ($r^L_m= r^{500}_m)_{m=1,\ldots, T}$ at the first $T=10 000$ time steps as an estimator for $\{r_m\}_{m=1,\ldots,T}$ (see a snapshot of $r^L_m$ at an arbitrary instance $m$ in Figure \ref{old_new_corr}). Using the same data assimilation output, we also estimate the correlations $\{r^K_m\}_{m=1,\ldots,T}$, subsampling $K$ out of $L=500$ members. In this numerical experiment, the ETKF is performed without any localization or inflation. Our choice of using ETKF is just to ensure that one can learn the map by fitting to correlation estimates from a different ensemble Kalman filtering method. In real applications, historical data or the output of a preferred cheap filter may be used as a proxy for the limiting correlation. The serial EnKF algorithm is then modified by applying the map $ r^K \leftarrow \mathcal{L}^T r^K $ at each assimilation cycle, which transforms the undersampled correlation matrix $r^K$ into a correlation matrix $\mathcal{L}^T r^K$ that is closer to $r$ (see the corresponding snapshots in Figure \ref{old_new_corr}).


\begin{figure}
\begin{center}
\subfigure[5 ensemble members.]{\includegraphics[width=10.5cm]{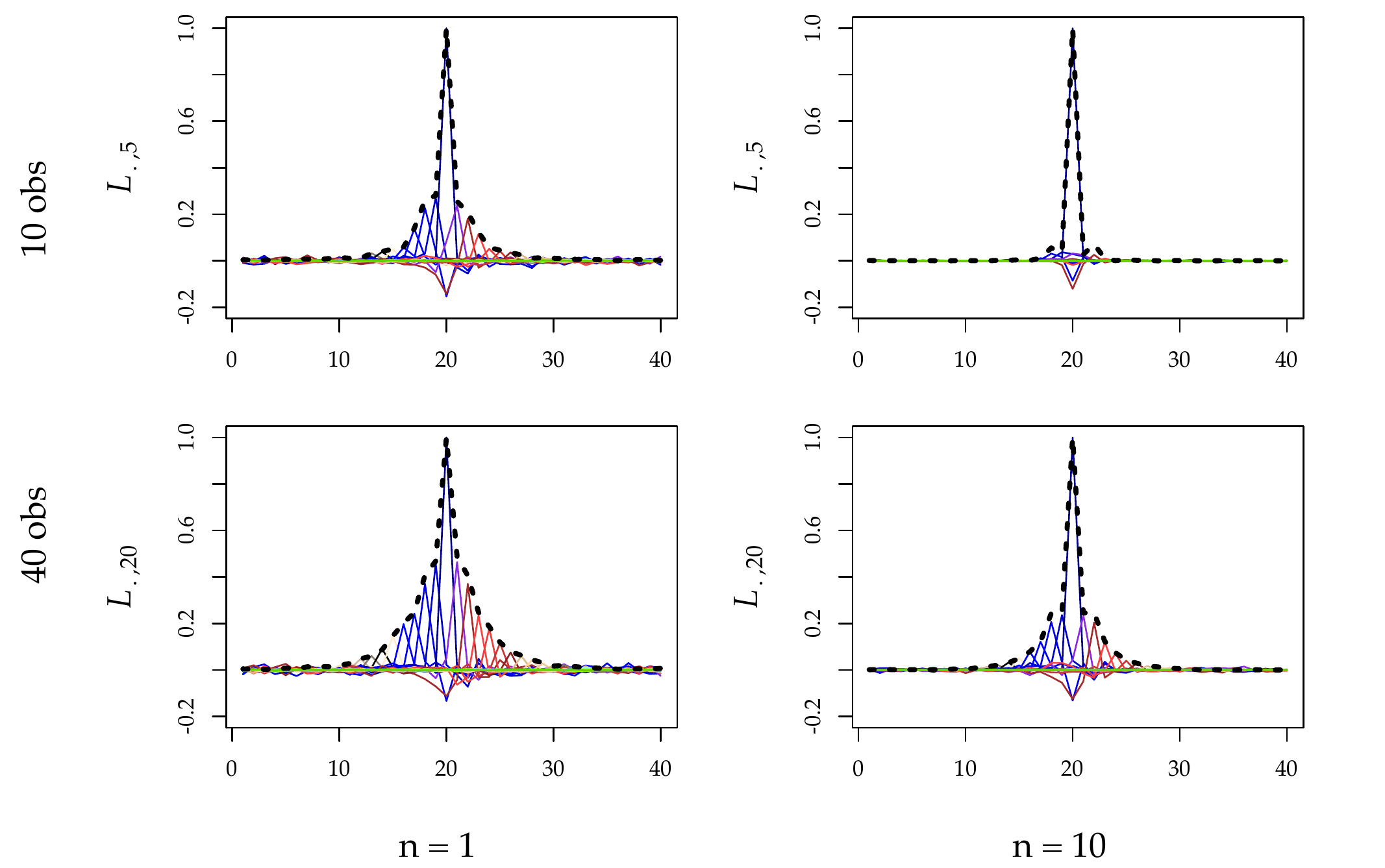}}  \subfigure[40 ensemble members.]{\includegraphics[width=10.5cm]{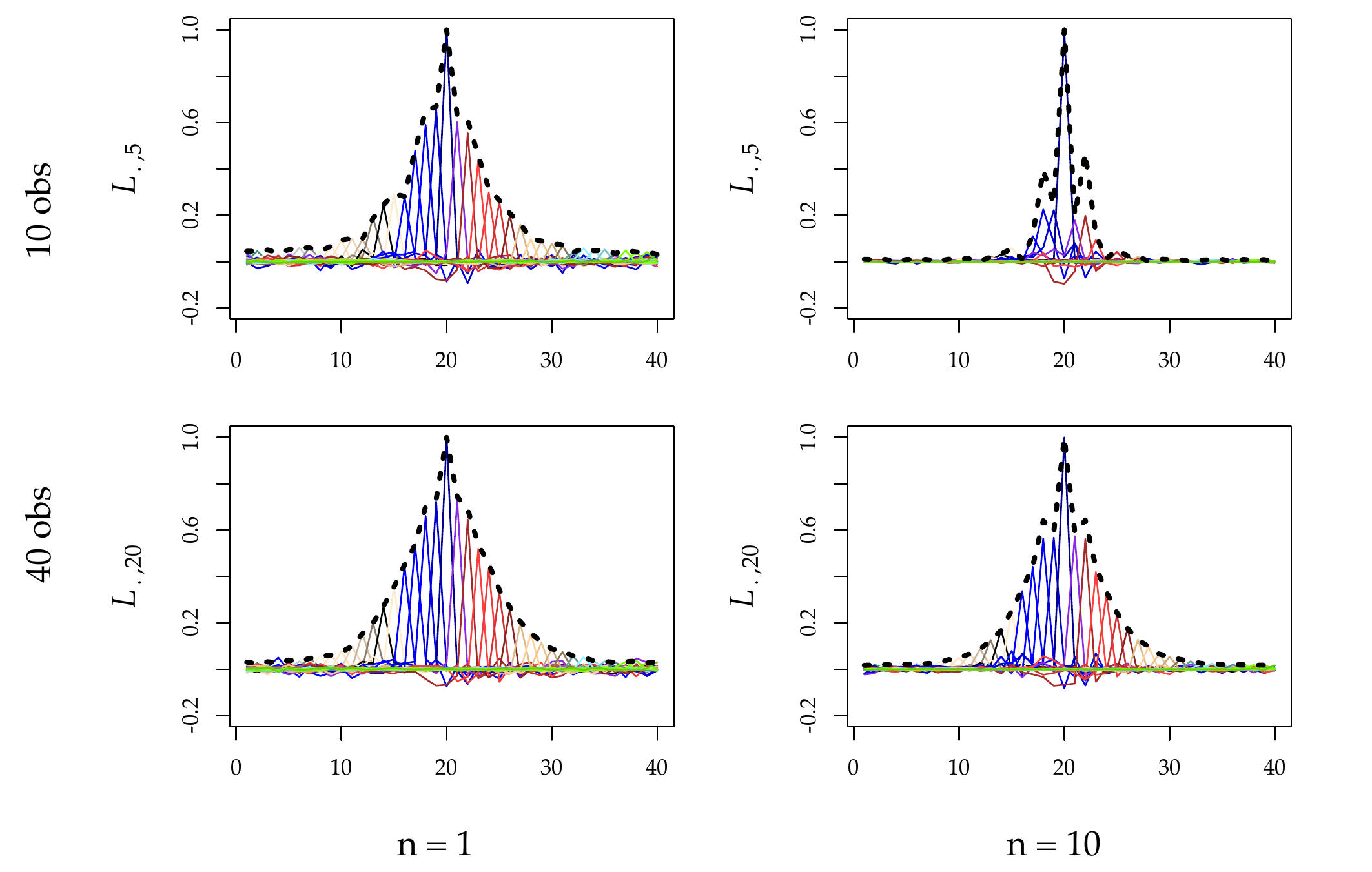}}\label{locmap_linear_direct_right}
\end{center}
\caption{The map $\mathcal{L}(q,i,j)$ for the middle observation $j=M/2$, shown as a function of the state variables $q=1,\dots,40$ for $M = $ 10 and 40 linear direct observations and for observation time steps $n=$ 1, 10. The map $\mathcal{L}$ is plotted here for every pair $(i,j)$ of state variable $i=1,\dots,40$ and the middle observation $j=M/2$ in the case of (a) 5 and (b) 40 ensemble members. The corresponding diagonal map $\mathcal{L}_d$ is overlaid in dashed black. }\label{locmap_linear_direct}
\end{figure}

\begin{figure}
\begin{center}
\hspace{-0.4cm}\subfigure[5 ensemble members.]{\includegraphics[width=15cm,height=5.5cm]{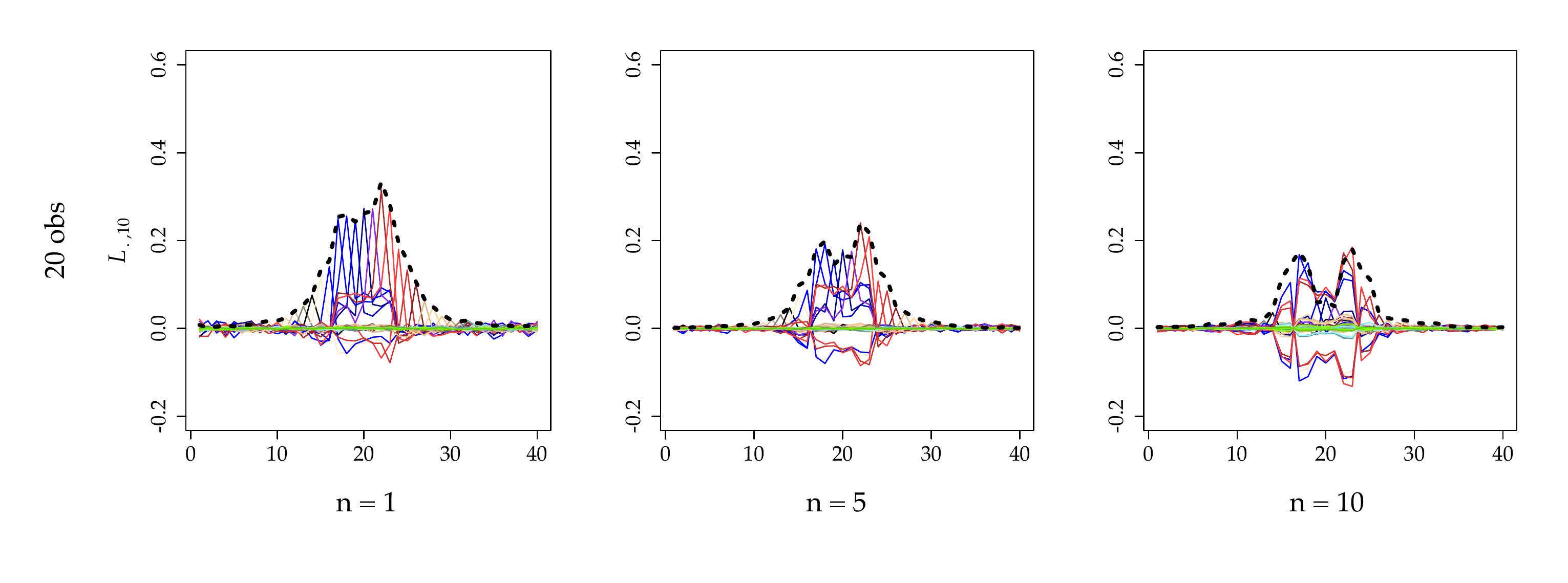}}

\hspace{-0.4cm}\subfigure[40 ensemble members.]{\includegraphics[width=15cm,height=5.5cm]{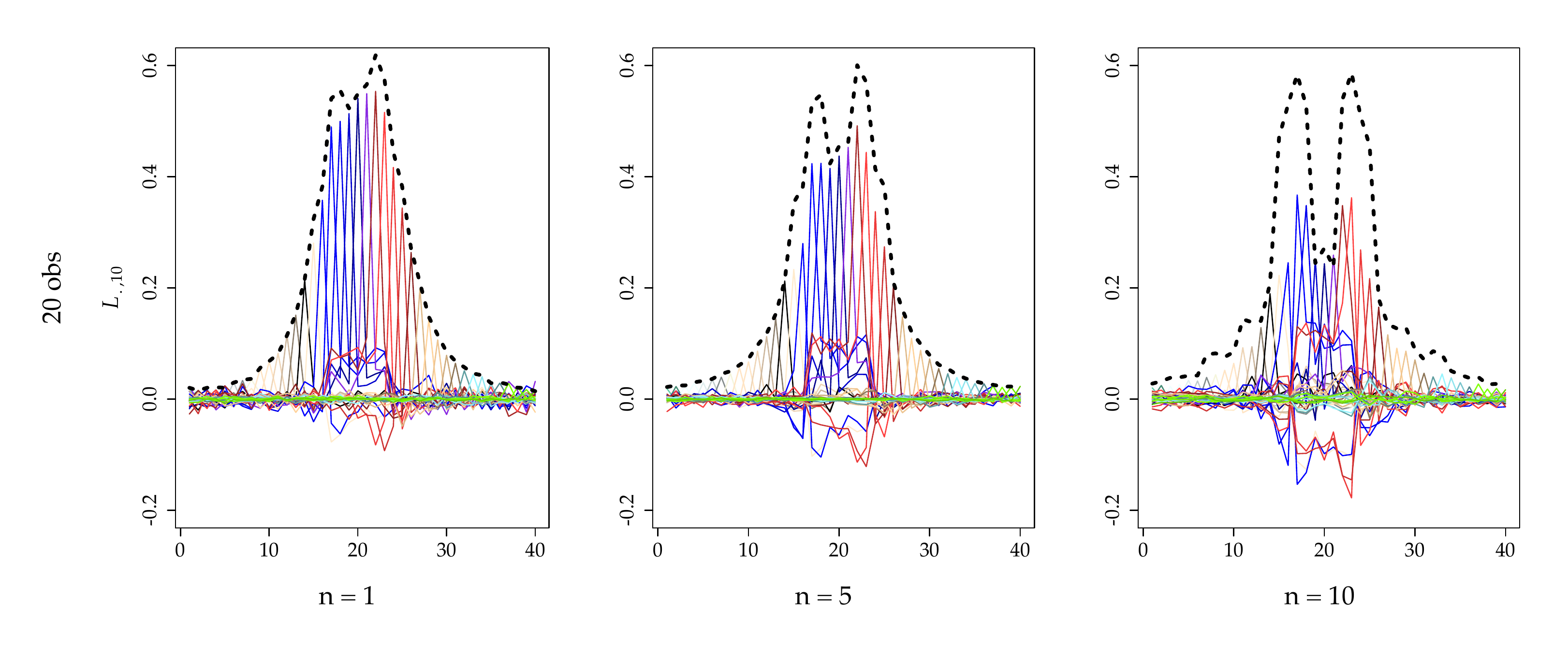}}
\end{center}
\caption{Same as in Figure \ref{locmap_linear_direct} but for $M=20$ indirect observations and observation time steps $n=1$, 5, and 10. }\label{locmap_linear_indirect}
\end{figure}

\textit{Linear Observations.}\quad The maps learned in the training phase for linear direct and indirect observations are shown in Figures \ref{locmap_linear_direct} and \ref{locmap_linear_indirect}, respectively. The maps $\mathcal{L}(\cdot,i,j)$ are plotted for every pair $(i,j)$ of state variables and observation at grid point $j=M/2$, where $M$ is the number of observations ($M= 10, 40$). These curves are compared with the diagonal maps which arise from a separate minimization procedure, fitting to a diagonal $\mathcal{L}(\cdot,\cdot,j) \in \mathbb{R}^{N\times N}$;  we refer to the resulting map as $\mathcal{L}_d$.
For direct observations, we note that both the map $\mathcal{L}$ and the diagonal map $\mathcal{L}_d$ peak to unity at the midpoint of the interval $[1,40]$, which occurs precisely for the pair $(20,M/2)$ of state and observation variables, and slowly decay away from the observed location. We also find that the map $\mathcal{L}_d$ can be viewed as an envelope for the packet of mapping curves $\mathcal{L}$. In this sense, the maps $\mathcal{L}$ and $\mathcal{L}_d$ capture the correlation structure between $\bs{x}$ and $H\bs{x}$, where $H$ is the direct observation operator. The mapping curves are also more compact or ``localized'' when the ensemble size used is small, or when the observation time step is large. Interestingly, the regression produces maps with a structure that resembles the localization functions that are commonly used in EnKF applications; the key point is that here the linear map is not chosen a priori but instead is learned from the data. Based on these observations, we will refer to  the maps $\mathcal{L}$ and $\mathcal{L}_d$ as \textit{localization map} and \textit{scalar localization function}, respectively. 

Although the structure of the localization mappings can be easily anticipated in the case of direct observations, it may be more challenging to design a localization function when the state variables are observed indirectly. As shown in Figure \ref{locmap_linear_indirect}, in the case of indirect observations, the localization maps  exhibit a nontrivial structure that may be difficult to guess a priori. Interestingly, the maps have larger weights when the ensemble size is 40, compared to an ensemble of size 5, especially for smaller $n$. These weights also tend to decrease in magnitude when the observation time step increases. This suggests that when the ensemble size $K$ is larger or when the observation time step $n$ is small, the localization map gives more weights to the sample correlation $r^K$ because it is a good estimator of $r$.  

\begin{figure}
\begin{center}
\subfigure[5 ensemble members.]{\includegraphics[width=6.2cm,height=5cm]{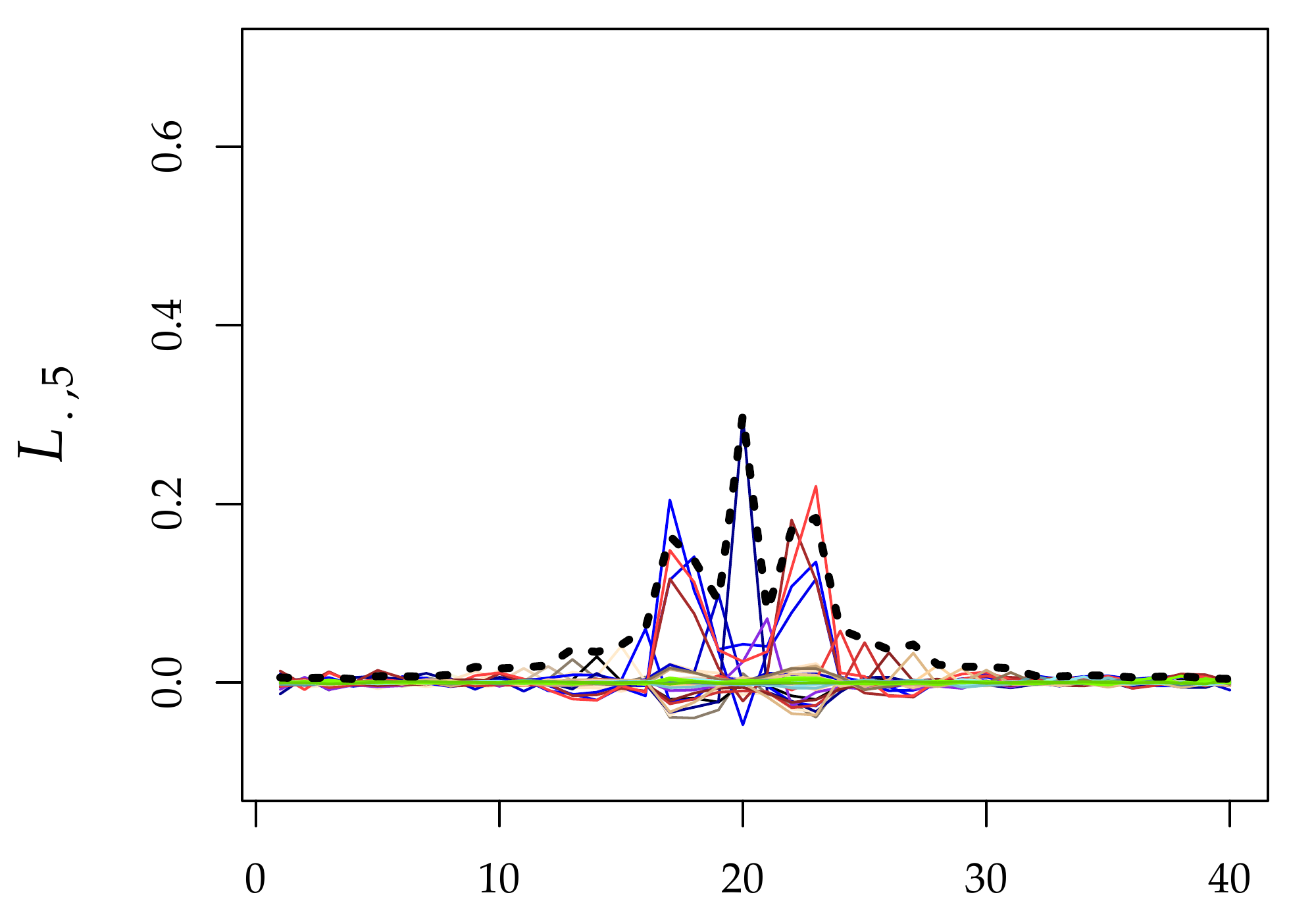}}
\hspace{1cm}
\hspace{-0.4cm}\subfigure[40 ensemble members.]{\includegraphics[width=6cm,height=5cm]{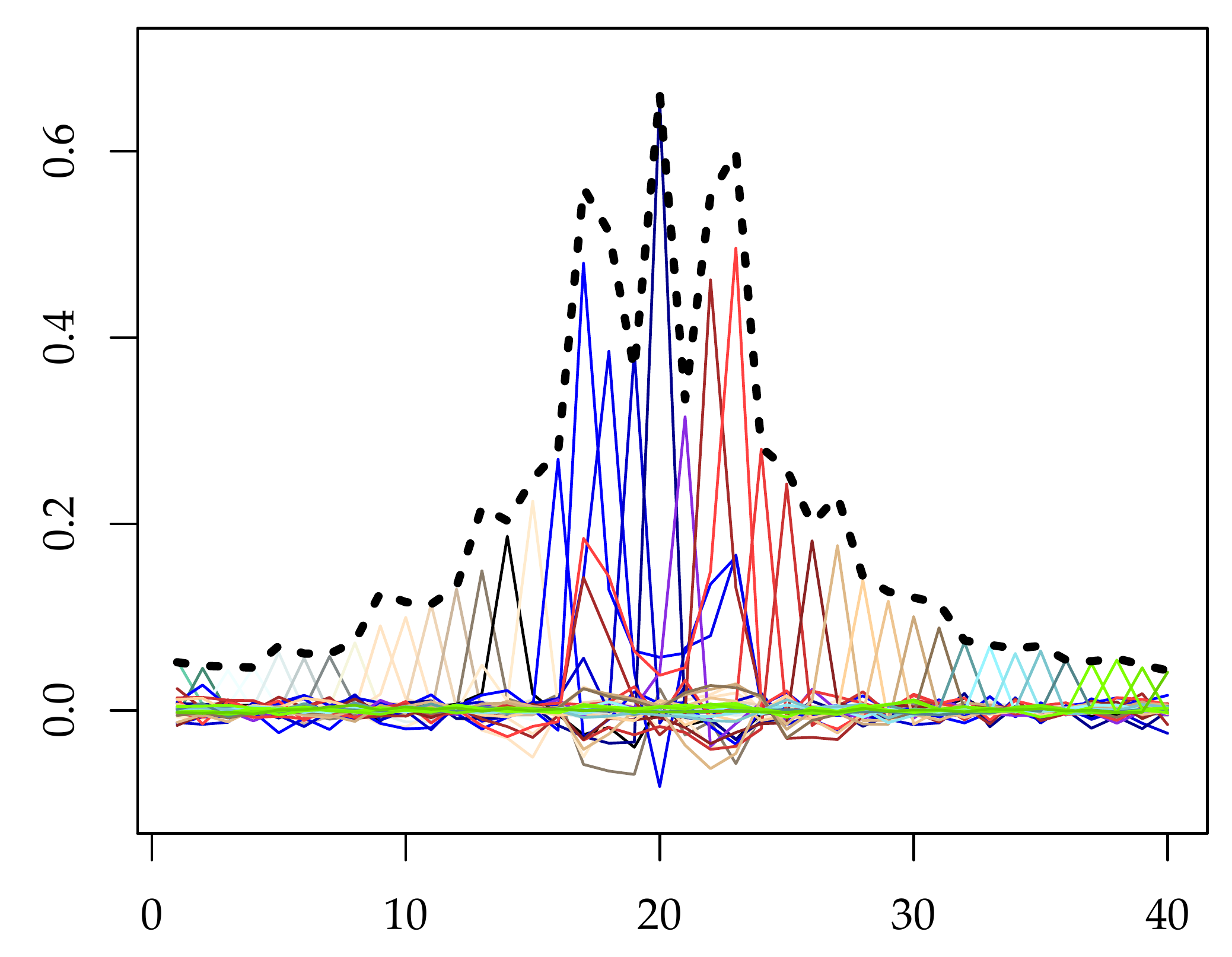}}

\end{center}
\caption{Same as in Figure \ref{locmap_linear_direct}, but for $M=10$ nonlinear indirect observations and observation time step $n=5$. }\label{locmap_nonlinear}
\end{figure}

\textit{Nonlinear Indirect Observations.}\quad The maps for 10 nonlinear indirect observations observed at every 5 model steps ($n=5$) are presented in Figure \ref{locmap_nonlinear}, for ensemble member of sizes 5 and 40. These maps show more complexity than the ones obtained in the linear indirect case, and reveal a particularly interesting trimodal envelope which is captured also in the diagonal maps $\mathcal{L}_d$. As in the linear indirect case, the maps yield much greater weights when the ensemble size is 40, compared to when the ensemble is only 5.

\begin{figure}
\begin{center}
\includegraphics[width=10cm]{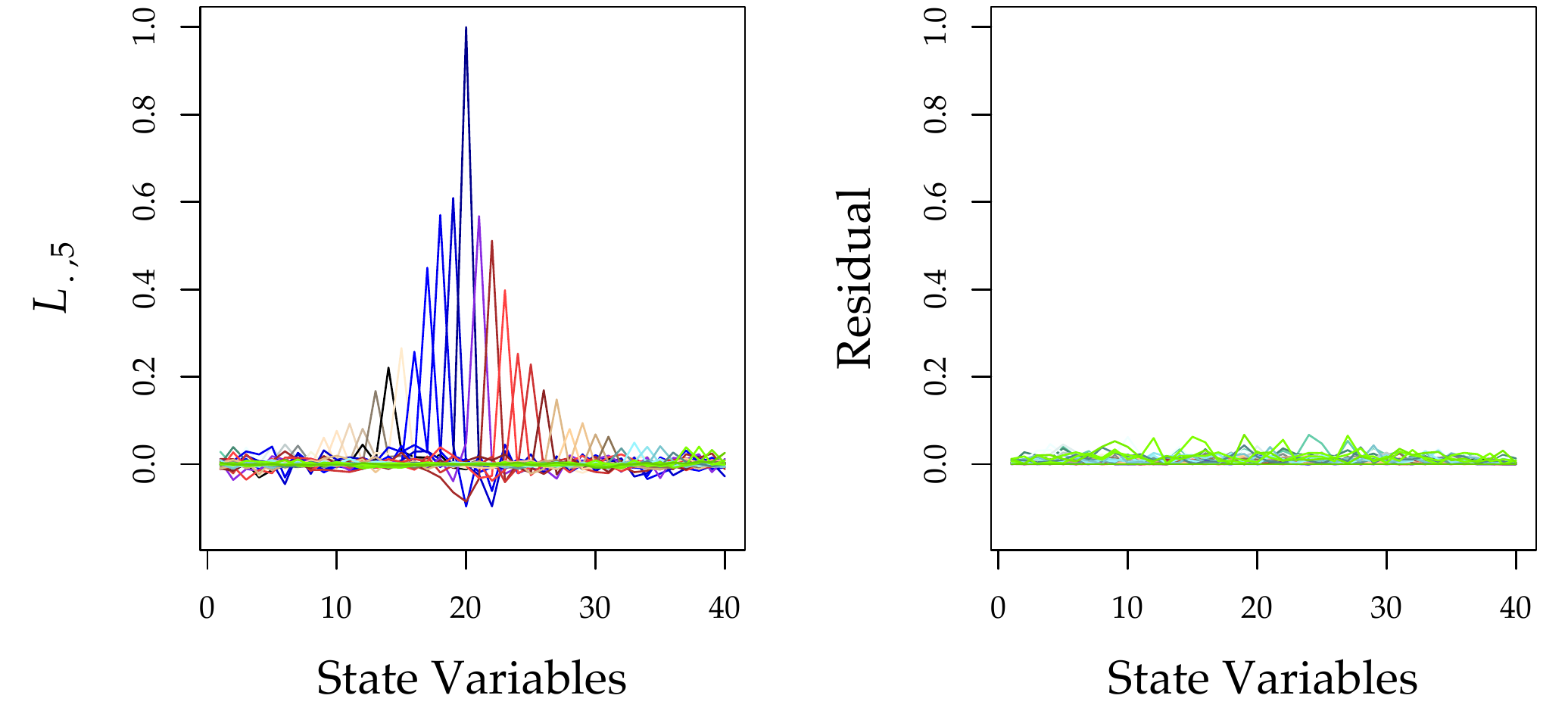}\end{center}
\caption{(\textit{Left}) Maps $\mathcal{L}$ obtained using the value $S=100$, in the case of 10 linear direct observations taken at every model step, and for ensemble size $40$. These maps should be compared with those of Figure 2b  ($M=10$ and $n=1$) obtained with $S=1$. (\textit{Right}) Residuals (pointwise difference) between the maps obtained with the values $S=1$ and $S=100$. }\label{residual}
\end{figure}

\textit{Robustness with respect to the parameter $S$.} 
It should be mentioned that, for the sake of computational efficiency, the localization maps shown in Figures \ref{locmap_linear_direct}, \ref{locmap_linear_indirect} and \ref{locmap_nonlinear} were all obtained for $S=1$
in the minimization problem (\ref{regression}), that is, when only one sample of the correlation $r^K$ is taken from the sampling distribution $p(r^K| r)$. However, as outlined in Section 3, as many as ${500 \choose K}$ samples can be obtained from randomly selecting $K$ members from the ETKF analysis ensemble $\{x^{a,k} \}_{k=1}^{500}$. 
In Figure \ref{residual}, we show the localization maps $\mathcal{L}$ obtained from the minimization procedure with $S=100$ correlation samples $\{r_s^K\}_{s=1}^{100}$, each one being computed by randomly selecting $K=40$ analysis members from the 500 available, in the case of 10 linear direct observations taken at every model step. The resulting maps are arguably similar, at least qualitatively, to the ones previously shown in Figure~\ref{locmap_linear_direct} for the value $S=1$, and the small residuals (pointwise difference) between the two sets of maps indicate that the method is not very sensitive to the parameter $S$. 
\vspace{0.2cm}

\subsection{Filtering results}

We now verify the filtering performance, assimilating the last 20 000 observations, beyond the first 10 000 data set that is used for training. 
The performance of the filter is measured using the following metrics:
\begin{itemize}
\item The root mean squared error (RMSE) of the analysis ensemble mean from the truth $\boldsymbol{x}^\dagger$, defined as
\begin{equation*}
RMSE = \sqrt{\dfrac{1}{N}\sum_{j=1}^{N} (\hat{\overline{x}}^{\,a}_j} - x^\dagger_{j})^2,
\end{equation*}
where $\hat{\bar{x}}_j^a$ and $x_j^\dagger$ denote the $j$th component of the true state and the analysis mean estimate, respectively.
\item The spread of the analysis from the ensemble mean, defined as
\begin{equation*}
spread = \sqrt{\dfrac{1}{N(K-1)}\sum_{j=1}^{N}\sum_{k=1}^{K} (\delta x_j^{k,a} )^2}   
\end{equation*}
where $\delta x_j^{k,a} =  x^{k,a}_j -  \hat{\overline{x}}^{\,a}_j$ is the deviation of the $k^{th}$ analysis ensemble member from its ensemble mean. 
\end{itemize}
Here we investigate the performance of the serial EnKF when the localization map $\mathcal{L}$ is applied in conjunction with a multiplicative covariance inflation. The covariance inflation factor $f$ is applied \citep{anderson:01} to the prior ensemble estimate of the state covariance as follows $(1+f)P^b$ in an attempt to reduce filter errors and avoid filter divergence.

\begin{figure}
\begin{center}
\subfigure[Localization map $\mathcal{L}$.]{\includegraphics[width=13cm]{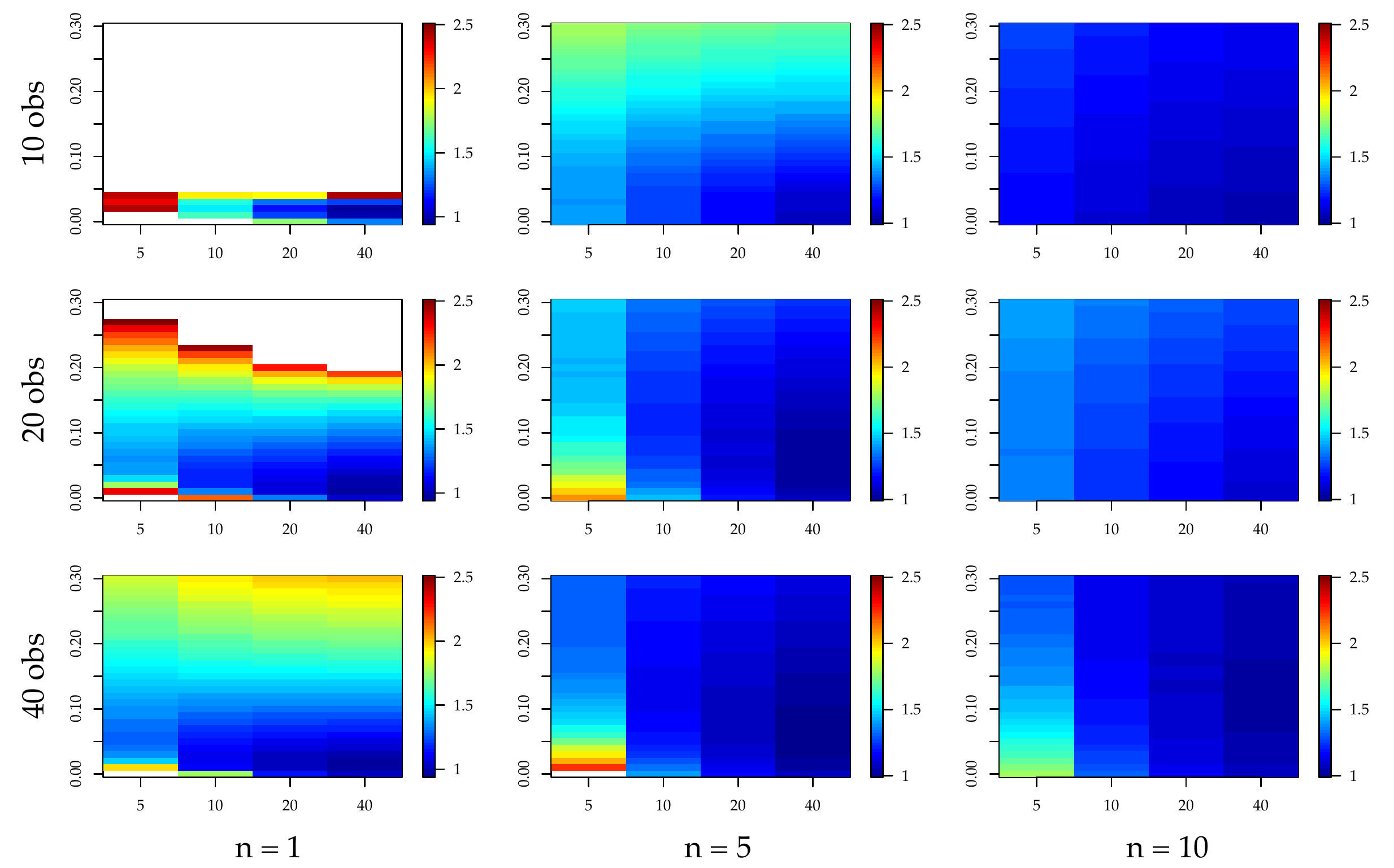}}\label{}
\subfigure[Scalar localization function $\mathcal{L}_d$.]{\includegraphics[width=13cm]{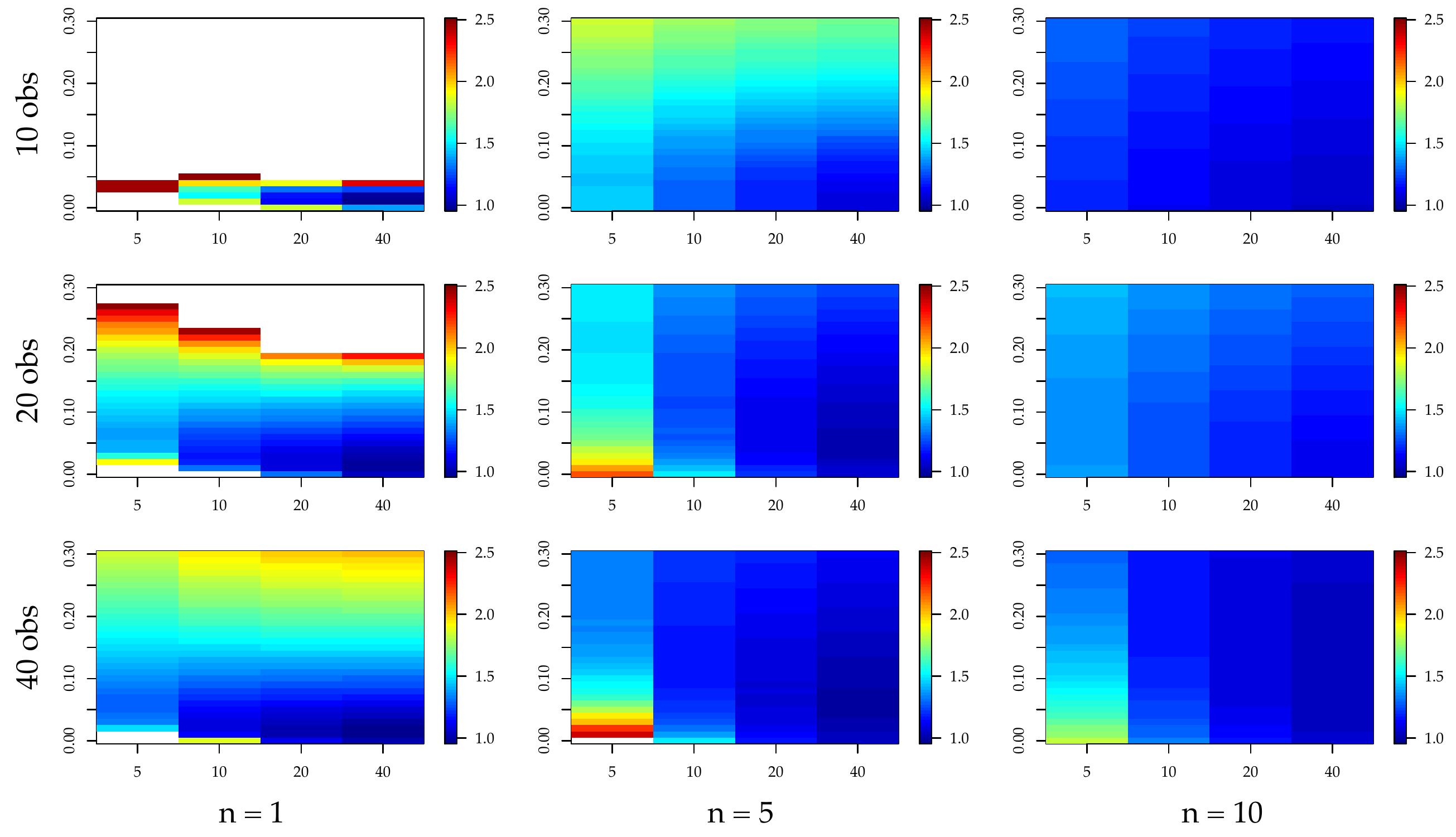}}\label{}
\end{center}
\caption{\textit{Linear direct observations.}  Array of time mean RMSE  colormaps (normalized by the RMSE of the ETKF using 500 ensemble members) for the serial EnKF using  (a)  the localization map $\mathcal{L}$ and (b) the scalar localization function $\mathcal{L}_d$, plotted against inflation factor and ensemble size $K=$ 5, 10, 20, 40. The 3 by 3 array's columns and rows correspond to observation time step $n$ and number $M$ of linear direct observations, respectively. White pixels correspond to RMSE values outside the scale shown.} \label{rmse_lineardirect}

\end{figure}

 \textit{Linear direct observations.} \, Figure  \ref{rmse_lineardirect} shows the time mean of the RMSE (averaged over the last 20 000 cycles) in the case of linear direct observations, as a function of ensemble size and covariance inflation factor. The  RMSE values reported here are normalized by the RMSE of the ETKF with 500 ensemble members, assimilating the same set of observations. Both the localization map $\mathcal{L}$ and scalar localization function $\mathcal{L}_d$ are tested on various number of observations ($M=10,20,40$) and observation time steps ($n=1,5,10$).

We see that compared to the benchmark ETKF solution, the serial EnKF filter using both the localization maps $\mathcal{L}$ and $\mathcal{L}_d$ performs well for large ensemble sizes, and its performance improves when the inflation factor is properly tuned. In the assimilation regime $M=10$, $n=10$, and $K=5$, the filter's error at its lowest is about 17\% more than the ETKF error, which is achieved when covariance inflation is completely switched off. Interestingly, if the observation time step is reduced to the model step ($n=1$), the accuracy of the serial EnKF analysis' solution degrades to more than twice the ETKF error (top left colormap), even when the inflation factor is well tuned. In fact, in the regime $M=10$ and $n=1$, the filter diverges when the inflation factor goes beyond a threshold value of approximately 0.1. However, in many assimilation regimes (e.g. $M=40$, $n=1$, and $K=40$), the serial EnKF performs just as well as the benchmark ETKF. Notice also that in these direct observation cases, the sensitivity of the filter's estimates to ensemble size and  inflation factor using the map $\mathcal{L}$ is not very different compared to that using the function $\mathcal{L}_d$ (compare colormaps (a) and (b) in Figure~\ref{rmse_lineardirect}).

\begin{figure}
\begin{center}
\subfigure[10 observations]{\includegraphics[width=11cm]{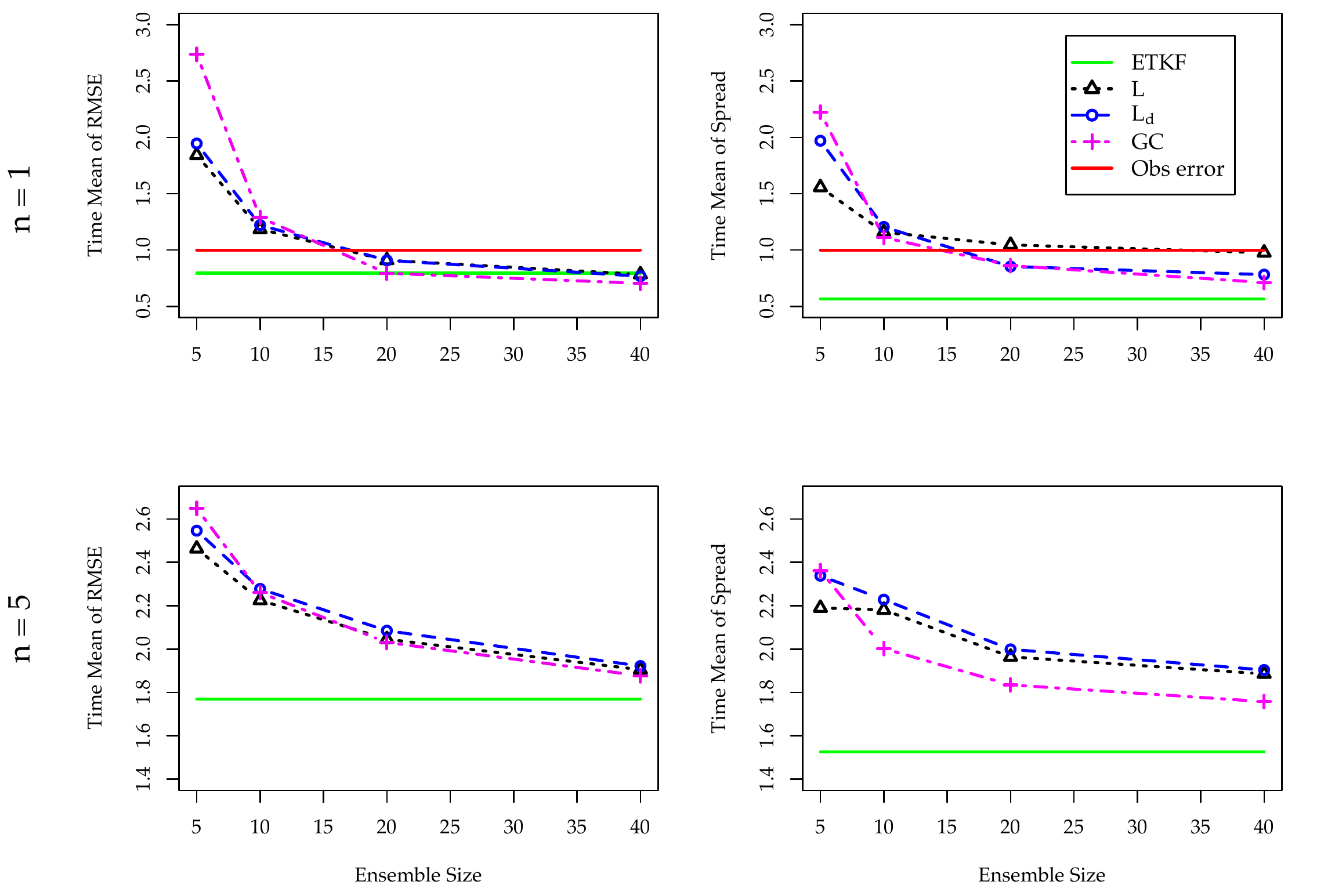}}
\subfigure[20 observations]{\includegraphics[width=11cm]{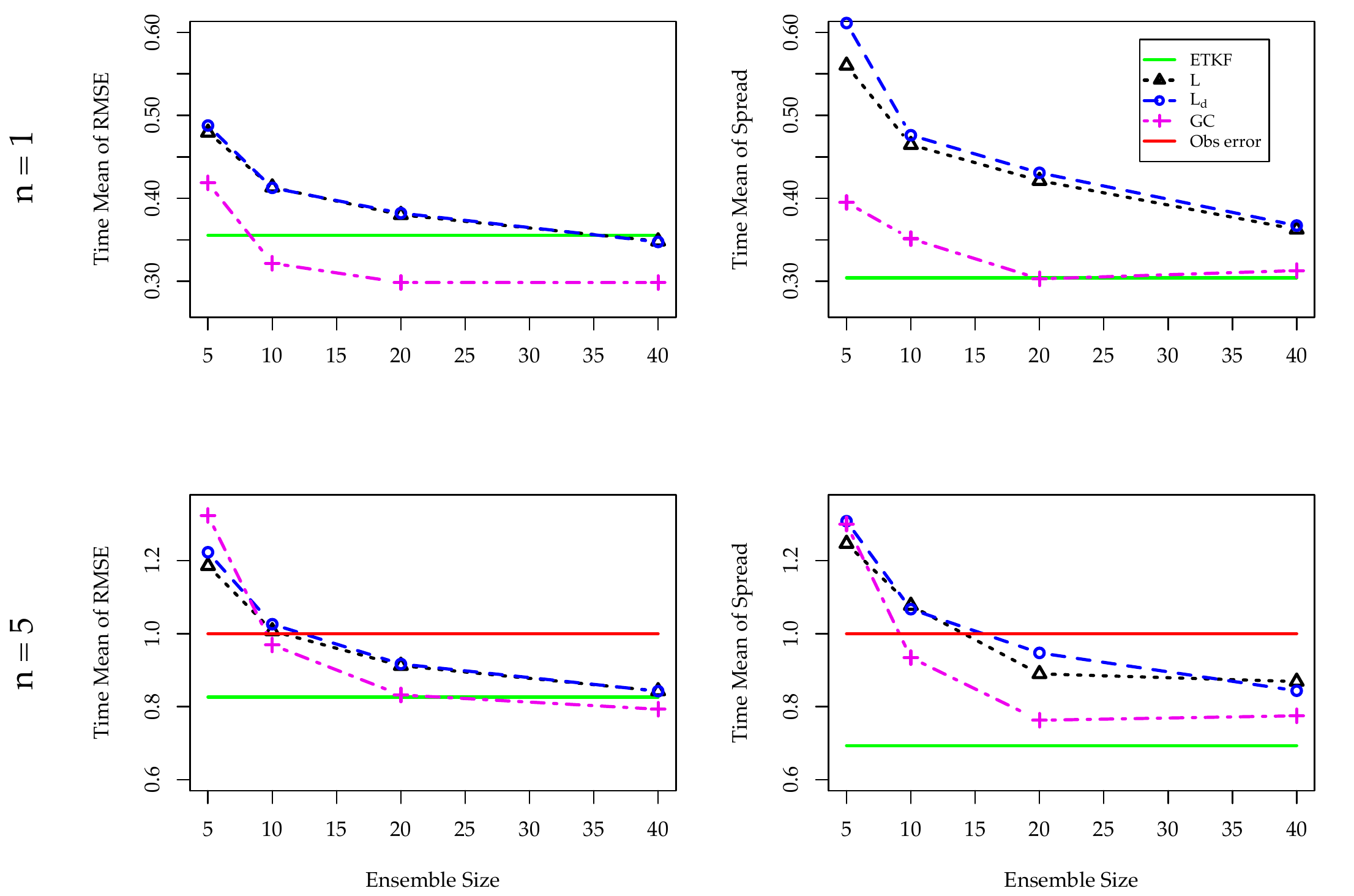}}
\end{center}
\caption{\textit{Linear direct observations.}  Time mean of RMSE and spread for the serial EnKF as a function of ensemble size using three  types of localization: the localization map $\mathcal{L}$, the scalar localization function $\mathcal{L}_d$, and the Gaspari-Cohn (GC) localization with an optimal half-width. All use the best tuned inflation values. These are compared with the benchmark ETKF solution using 500 ensemble members. Also shown is the observation error variance.} \label{directobs_5curves}
\end{figure}

In Figure \ref{directobs_5curves}, the time mean of RMSE for the best tuned inflation values are shown as a function of ensemble size, in the case of 10 and 20 direct linear observations, and observation time steps $n=1$ and $n=5$. In the same figure, we also show the corresponding time mean of spread. The figure compares three types of localization methods for the serial EnKF: the localization map $\mathcal{L}$, the scalar localization function $\mathcal{L}_d$, and the standard Gaspari-Cohn localization with an optimal half-width that is tuned in the training phase over a range of 1 to 10 grid points in each assimilation case (each case consisting of a specific choice of $K$, $n$, $M$, and inflation factor). We note that in most cases, the filter's performance using the localization maps $\mathcal{L}$ and $\mathcal{L}_d$ is comparable with that using the GC localization function. In fact, when only 5 ensemble members are used, the localization maps $\mathcal{L}$ and $\mathcal{L}_d$ outperform GC in all instances shown except for the case $M=20$, $n=1$. In this later case, and for all larger ensemble sizes, GC still gives a more accurate filter solution. Nevertheless, these are still significant results, considering that the maps $\mathcal{L}$ and $\mathcal{L}_d$ are obtained offline without any tuning.

We found that the optimal GC full widths for the case of $n=1, M=10$ are $6, 16, 12, 16$ for $K=5, 10, 20, 40$, respectively. Comparing these optimal GC widths with the support of the localization maps in  Figure~\ref{locmap_linear_direct} for the cases of $K=5$ and $K=40$, one can see that the map is negligible outside the support of the GC function.

In our cases of interest, i.e. for the localization maps $\mathcal{L}$ and $\mathcal{L}_d$, the spread is close to the RMSE in most cases (but not always): the ensemble spread is large (more uncertainty about the estimator) when the ensemble mean lies far from the truth, and the spread is small (more confidence about the estimator) when the ensemble mean lies close to the truth. 

\begin{figure}
\begin{center}
\includegraphics[width=10cm]{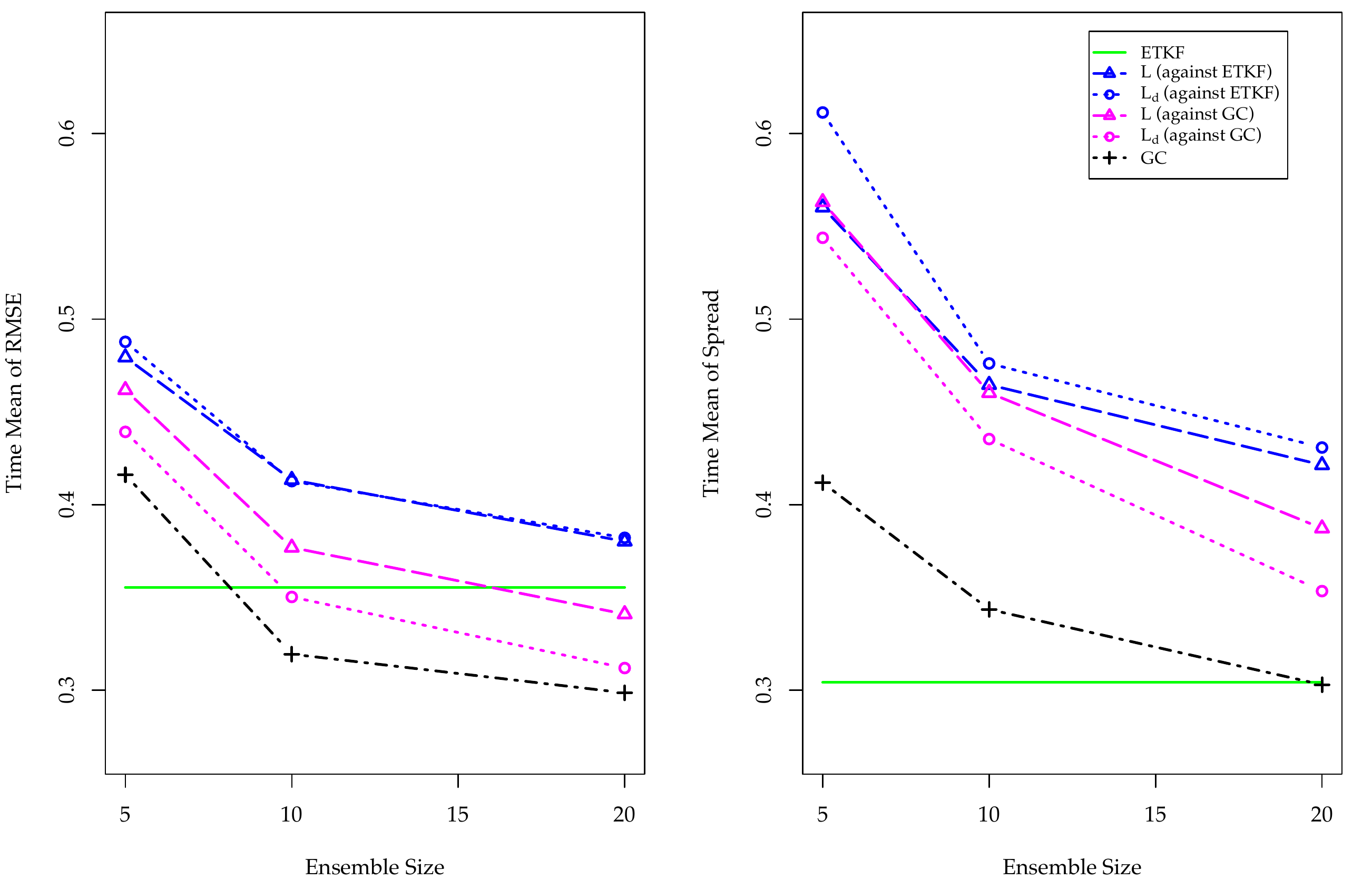}\end{center}
\caption{Time mean of RMSE and spread as a function of ensemble size for the serial EnKF with localization maps $\mathcal{L}$ and $\mathcal{L}_d$, obtained from two different choices of regressors: 1) the ETKF with 500 ensemble members, and 2) the serial EnKF with GC localization with optimal half-width. All use the best tuned inflation values. The results are for 20 linear direct observations taken at every model step. }\label{M20_againstGC}
\end{figure}

From the numerical results of Figure \ref{directobs_5curves}, notice that the serial EnKF with a well-tuned GC localization beats the ``benchmark'' ETKF in many cases when the ensemble size is large (closer to 40). At this point, one could ask whether a different choice of regressor for the localization maps $\mathcal{L}$ and  $\mathcal{L}_d$, such as the one obtained from the output of the serial EnKF with a well tuned GC localization, could produce an even more accurate filter solution. In Figure \ref{M20_againstGC}, we show the results using both the ETKF and the serial EnKF with GC localization as regressors, in the case of 20 linear direct observations taken at every model step. In the case of the GC regressor, the maps are trained on the correlation $r$ obtained from 40 ensemble members of the serial EnKF with best tuned GC localization half-width parameter. The regressant $r^K$ is then computed by selecting 5, 10, or 20 members out of these 40 ensemble members. Interestingly, the filter's solution is more accurate using this new (GC) regressor, and the gain in accuracy is even more notable in the case of the scalar localization function $\mathcal{L}_d$.  This is in fact a meaningful test, since it demonstrates how different regressors, such as the one obtained from the output of a cheap preferred filter, may be used without undermining the filtering skills of the localization map method. 
 
\begin{figure}
\begin{center}
\subfigure[Localization map $\mathcal{L}$.]{\includegraphics[width=13cm]{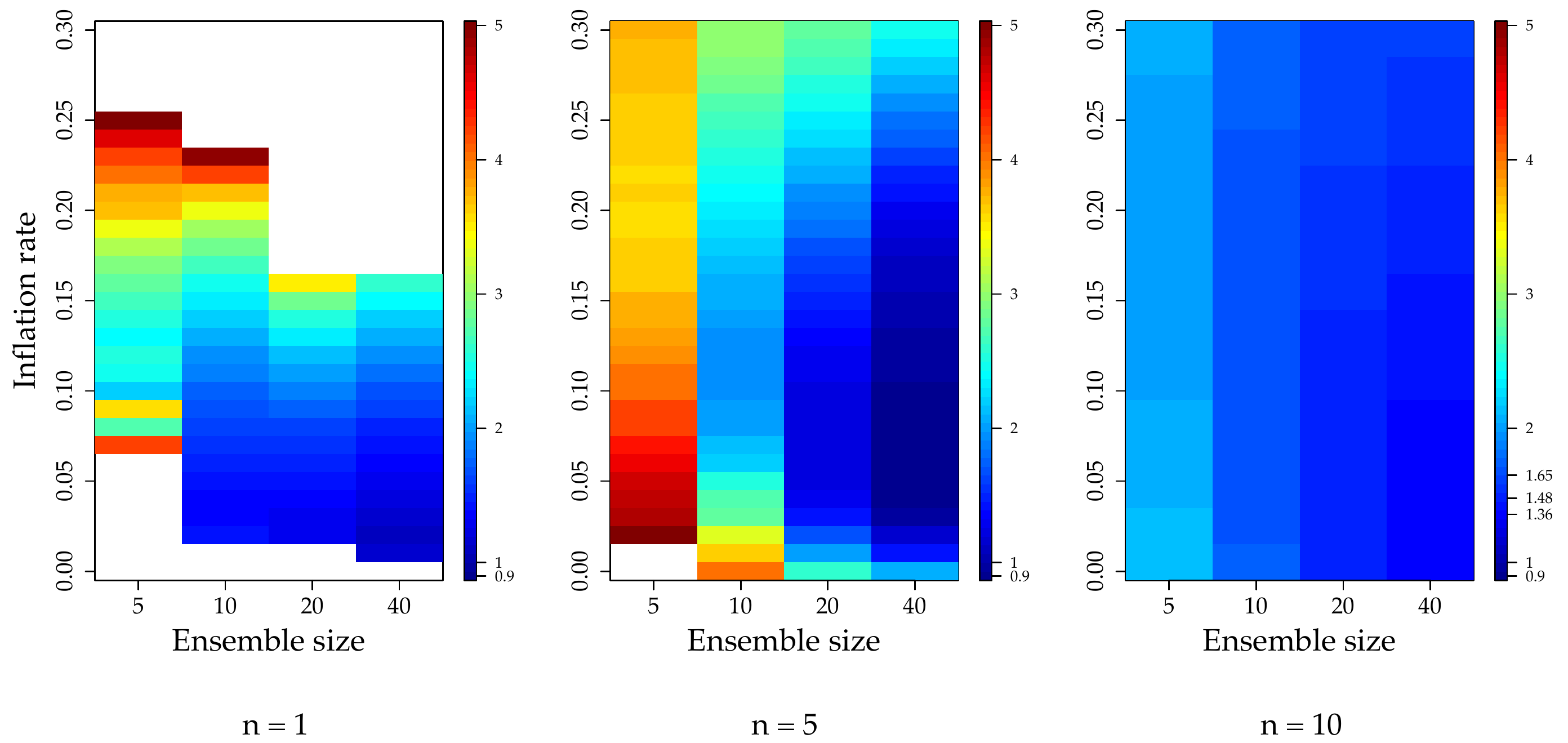}} \label{}\subfigure[Scalar localization function $\mathcal{L}_d$.]{\includegraphics[width=13cm]{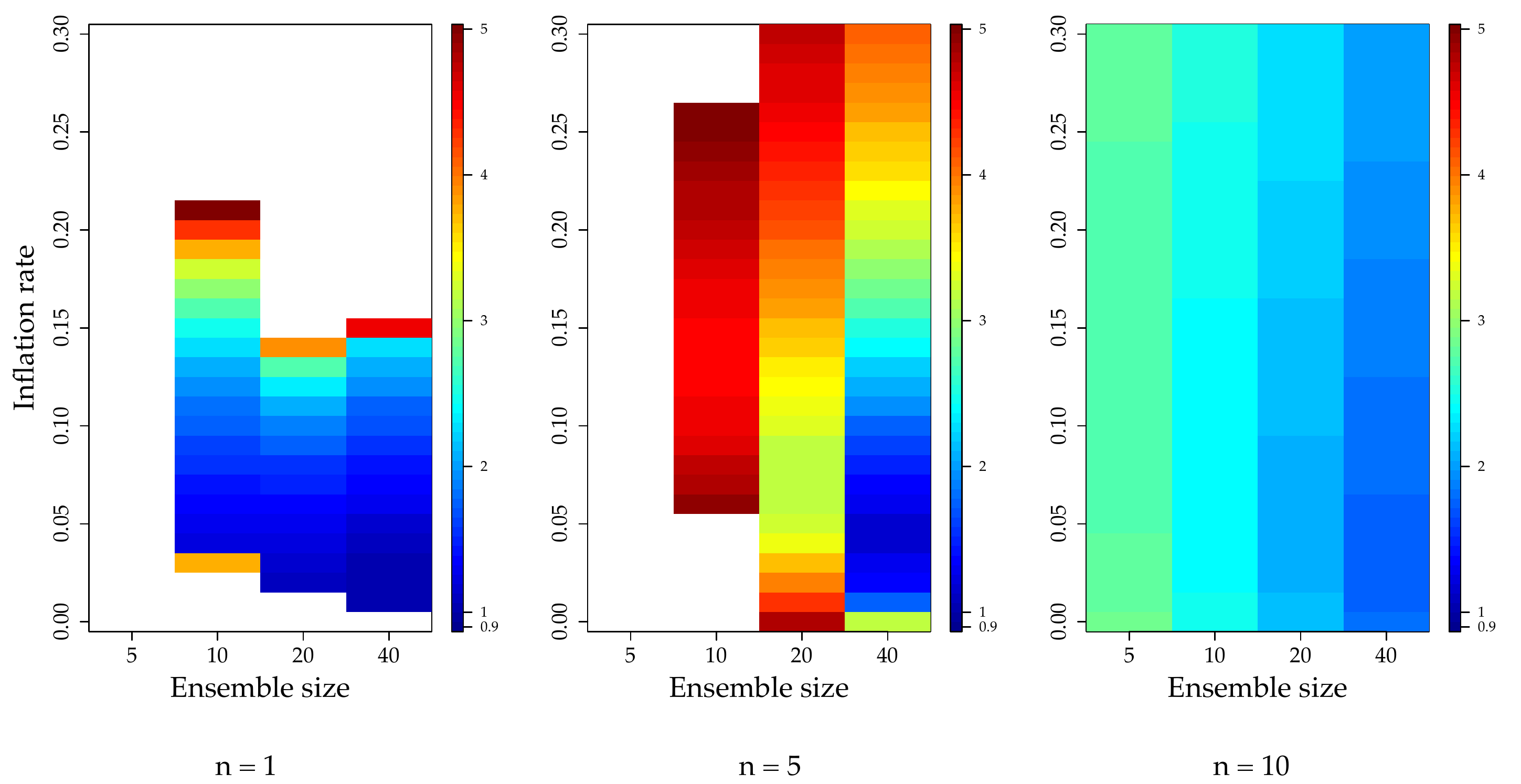}}
\end{center}
\caption{\textit{Linear indirect observations.}  Same as in Figure \ref{rmse_lineardirect}, but for 20 linear indirect observations and with a different scale. }\label{RMSE_linearindirect}
\end{figure}

\begin{figure}
\begin{center}
\includegraphics[width=12cm]{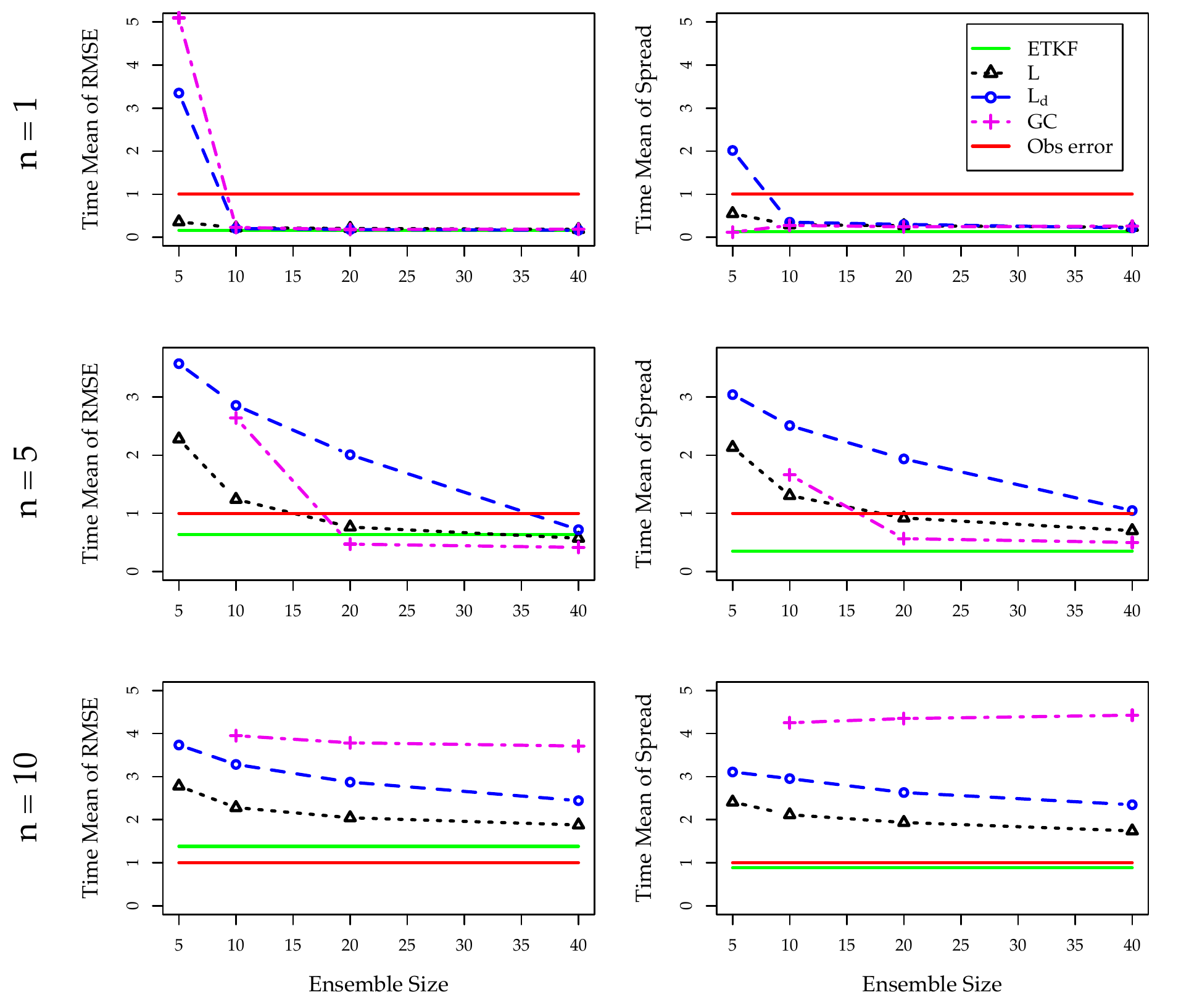}\end{center}
\caption{\textit{Linear indirect observations.}  Same as in Figure \ref{directobs_5curves}, but for 20 linear indirect observations, and observation time steps $n=1$, 5, and 10.  }\label{5curves_linearindirect}
\end{figure}

\textit{Linear indirect observations.} \, We next show in Figures \ref{RMSE_linearindirect} and \ref{5curves_linearindirect} the results when 20 linear indirect observations are taken at observation time steps $n=1,5$ and 10. It is worth noting that when observations are taken at every model step ($n=1$) and only 5 ensemble members are used, the performance of the map $\mathcal{L}$ is remarkable compared to that of the GC localization: The RMSE for $\mathcal{L}$ and GC are respectively 0.3602 and 5.0970, while the benchmark ETKF is at 0.1626  (see Figure \ref{5curves_linearindirect}). In this difficult regime, even the diagonal map $\mathcal{L}_d$ outperforms GC with an RMSE of 3.3498. If one uses 10 members or more,  all three localization methods perform qualitatively just as well as ETKF. For instance, in the case of $K=10$, the RMSE for $\mathcal{L}_d$, $\mathcal{L}$ and GC are respectively 0.2033, 0.2182, and 0.2276, compared to 0.1626 for ETKF. 
 
 When observations are taken at every  5 model steps $(n=5)$ and the ensemble size is 5 or 10, the map $\mathcal{L}$ offers a clear advantage over the other two  methods. In fact when only 5 members are used, the filter with a GC localization produces a divergent solution, while the RMSE for the map  $\mathcal{L}$ is 2.2793, nearly four times the ETKF error (0.6369). In comparison, the diagonal map $\mathcal{L}_d$ gives an error of 3.5728. However, for larger ensemble sizes, GC is outperforming the map $\mathcal{L}$, albeit not by much: for $K=20$, the RMSE for GC, $\mathcal{L}$, and $\mathcal{L}_d$ are respectively 0.4708, 0.7652, and 2.0077. 
For $K=40$ the errors are, in the same order, 0.4145, 0.5722, and 0.7191. 
 
When observations are taken at every $n=10$ model steps, the localization maps $\mathcal{L}$ and $\mathcal{L}_d$ both outperform the GC localization (see Figure \ref{5curves_linearindirect}). In particular, the well-tuned GC estimates diverge with RMSEs above the climatological standard deviation (approximately $ 3.6$), whereas the estimates from $\mathcal{L}$ are much more accurate (or below the climatology). Also, as shown in Figure \ref{RMSE_linearindirect}, the filter with either a localization map $\mathcal{L}$ or $\mathcal{L}_d$ produces a solution that is quite robust with respect to the inflation rate. In the case of the map $\mathcal{L}$ and when 5 ensemble members are used, 
the error produced is about twice that obtained with the ETKF. That ratio goes down to about 1.65 when $K=10$, 1.48 when $K=20$, and 1.36 when $K=40$. When the filter is used with the diagonal map $\mathcal{L}_d$, the  error climbs to nearly 3 times the error of the benchmark ETKF in the case of 5 ensemble members.

\begin{figure}
\begin{center}
\includegraphics[width=12cm ]{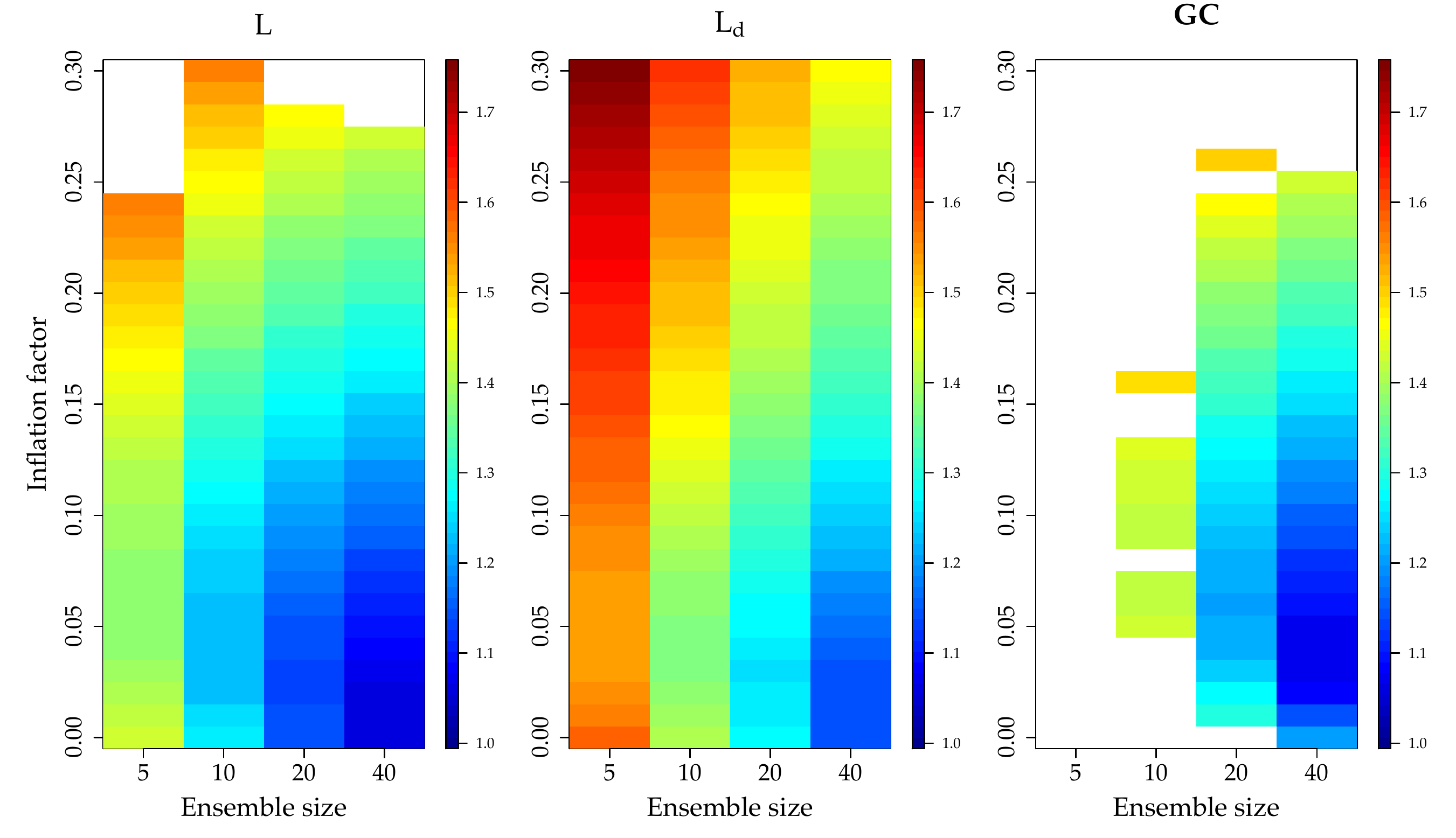}
\end{center}
\caption{Time mean RMSE for the serial EnKF in the case of 10 nonlinear indirect observations and observation time step $n=5$, using three different localization methods:  the localization maps $\mathcal{L}$ and $\mathcal{L}_d$, and the Gaspari-Cohn (GC) localization function with optimal half-width. The RMSE is normalized by the RMSE of the ETKF using 500 ensemble members, and plotted against inflation factor and ensemble size. White pixels correspond to RMSE values outside the scale shown.}\label{nonlinear2_llsenkf_rmse}
\end{figure}

\begin{figure}
\begin{center}
\includegraphics[width=11cm ]{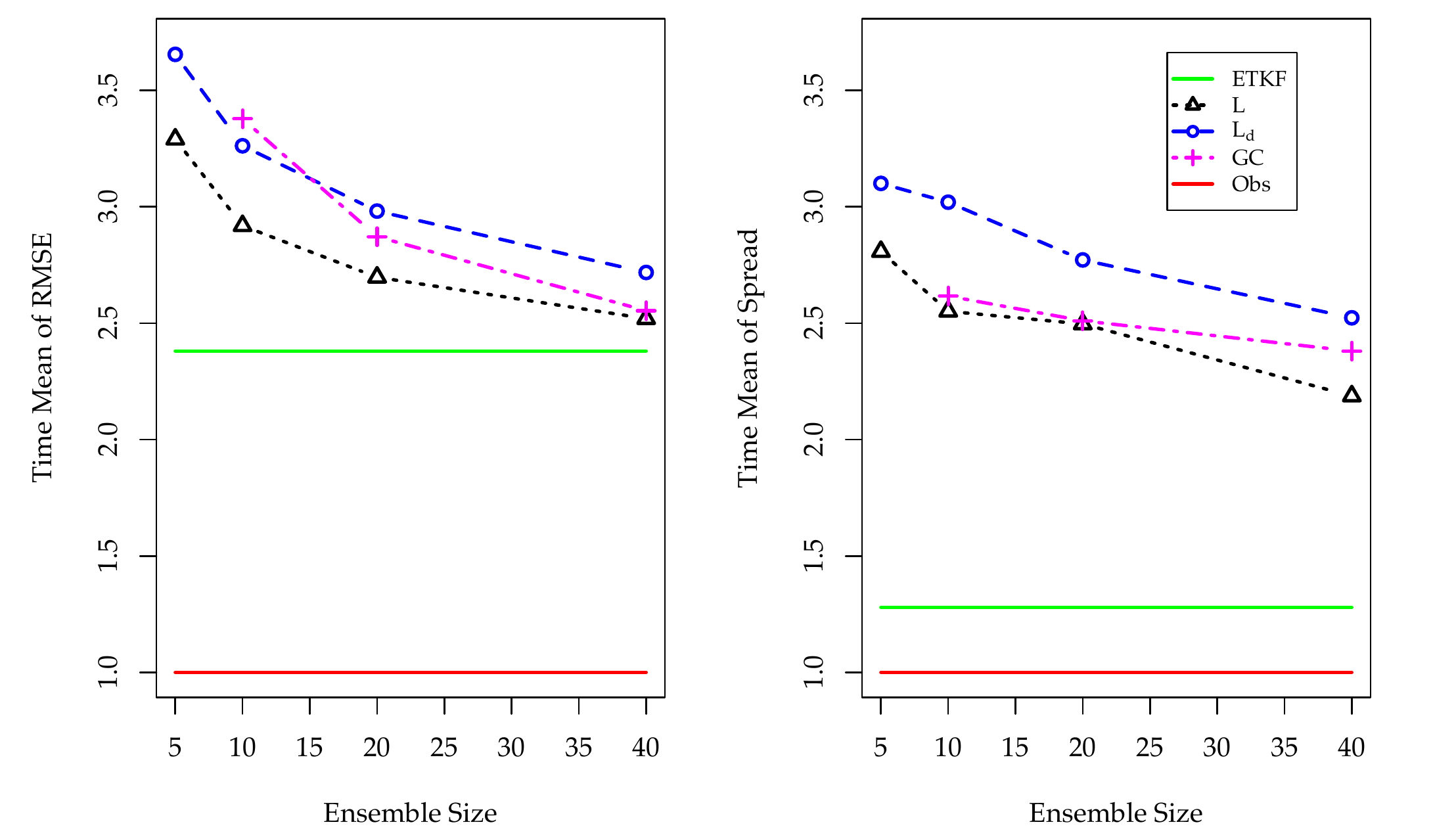}
\end{center}
\caption{\textit{Nonlinear Indirect Observations.}  Same as in Figure \ref{directobs_5curves}, but for 10 nonlinear indirect observations and observation time step $n=5$.  }\label{nonlinear_5curves_M10}
\end{figure}

 \textit{Nonlinear indirect observations.} \, In Figures \ref{nonlinear2_llsenkf_rmse} and \ref{nonlinear_5curves_M10}, we compare the results for 10 nonlinear indirect observations as defined in \eqref{nlinobs} with observation time step $n=5$. In this experiment, the localization map $\mathcal{L}$ outperforms both $\mathcal{L}_d$ and the well tuned GC localization for all ensemble sizes used. In Figure~\ref{nonlinear2_llsenkf_rmse}, we also include the sensitivity of the GC estimates with respect to inflation factor and ensemble size in the verification phase. As in all the experiments above, we should mention again that the GC localization radius is obtained in the training phase. Notice that 
when only 5 ensemble members are used ($K=5$), the filter's solution with a GC localization diverges catastrophically (solution blows up), while the filter's estimate gives an RMSE of 3.29 when a localization map $\mathcal{L}$ is used. As the ensemble size increases, the filter estimates become closer to that of the benchmark ETKF. Even for $K=10$, the GC filter's estimate is very sensitive with respect to variance inflation.

\section{Summary and conclusions}\label{section5}
Ensemble Kalman filters provide the exact filter solution for linear models with Gaussian observational error in the limit of large ensemble sizes \citep{mandel:12}. When the ensemble size is small, sampling errors tend to induce spurious correlations between observations and state variables. Localization methods address this issue by limiting the impact of an observation on a state variable, generally by specifying a radius of influence based on the physical distance between the two. However, most localizations require extensive tuning, which becomes impractical in large applications with multiscale spatial correlations \citep{dong:11,ki:14,flowerdew:15}. Moreover, when the observations are not correlated to the state variables by spatial distance, as may be the case in nonlinear models, the design of such localizations can be very difficult. 

In this work we propose a data-driven localization method for the serial EnKF which assumes no functional dependence on physical distance or other imposed restrictions, and more importantly, it does not require any tuning. The localization is in fact ``learned'' from any time series of correlations -- the regressor -- that  constitutes a good estimate of the filter limiting correlation, and is given as the solution of an appropriate linear regression problem. The localization takes indeed the form of a linear map $\mathcal{L}$ that transforms, at each EnKF cycle, the poorly estimated sample correlation into an improved correlation. We numerically tested the proposed methodology in an OSSE using the 40-variable Lorenz-96 model, with three different observation configurations: linear direct, linear indirect, and nonlinear indirect observations. We use as regressors the assimilation output of an ETKF with 500 ensemble members, and of a serial EnKF with a well-tuned GC localization with 40 ensemble members.  

For linear direct observations, the linear maps $\mathcal{L}$ and their diagonal counterparts $\mathcal{L}_d$ have a well defined structure that captures the intrinsic correlations of the filter solution: they peak to unity at every observed state variable, and slowly decay away from that location. This structure is reminiscent of the localization functions that are commonly used in EnKF localization. For the indirect and nonlinear indirect observations, we obtain localization maps with bimodal and trimodal envelopes, respectively. From our experiments, we found that when the true correlation is local (as it is in the linear direct observation example) our method is not advantageous compared to the well-tuned GC. On the other hand, when the correlation is nonlocal (e.g. in the linear indirect and nonlinear indirect observation examples), our method consistently performs well, beating the well-tuned GC, since the localization map contains such nontrivial correlation information.

In the case of the ETKF regressor, we found that the linear map $\mathcal{L}$ is more accurate and robust compared to both  $\mathcal{L}_d$ and the well-tuned GC localization method. In the case of the GC regressor, we found that the accuracy of the filter solutions with the localization maps $\mathcal{L}$ and $\mathcal{L}_d$ is improved compared to that obtained with the ETKF regressor. This test demonstrates that 1) the performance of the filter using the localization maps can be improved by using an even better correlation output as regressor and 2) the output of a cheap filter with small ensemble sizes can be used as a regressor when large ensemble data assimilation are too costly to be simulated, especially in large-dimensional applications. Therefore, the learning methodology proposed in this paper can be adopted using any available data assimilation correlation statistics and the results depend on the quality of the training data, like in any data-driven method.  Of course, the resulting map is also subject to noise but one can increase the number of data by fitting to the correlations with similar structure at different locations in the training phase. Finally, we should point out that one can use the nontrivial structure obtained from the training procedure here as a guideline to design a more appropriate parametric localization function.

While the localization map methodology was developed for the serial EnKF by correcting the correlations between observations and grid points based on learned data, we suspect that the same approach can be used on other variants of ensemble Kalman filters to improve the correlation statistical estimation. Also, the experiments reported in this paper are based on the perfect model assumption, while the proposed fitting strategy is not limited to OSSE. Thus it will be interesting to study the robustness of the proposed method in the presence of model errors.  We plan to address these issues and implement this method on higher dimensional applications in our future research.


\acknowledgments
This research is partially supported by ONR MURI
Grant No. N00014-12-1-0912, ONR Grant No. N00014-13-1-
0797, and NSF Grant No. DMS-1317919.


\bibliographystyle{ametsoc2014}
\bibliography{biblio}



\end{document}